\documentclass[11pt]{article}
\usepackage[latin1]{inputenc}
\usepackage[english]{babel}
\usepackage{amssymb,amsfonts, amsmath, color}
\usepackage[margin=2.7cm]{geometry}
\usepackage{xcolor}
\usepackage{graphicx, color, enumerate}
\usepackage[latin1]{inputenc}
\usepackage[active]{srcltx}
\usepackage{tikz}
\usepackage{pgf}
\usepackage{etex}
\usepackage{verbatim}
\usepackage{tikz-3dplot}
\usepackage{pgfkeys}
\usepackage{fp}

\newtheorem{theorem}{Theorem}[section]
\newtheorem{proposition}[theorem]{Proposition}
\newtheorem{corollary}[theorem]{Corollary}
\newtheorem{lemma}[theorem]{Lemma}
\newtheorem{remark}[theorem]{Remark}
\newtheorem{definition}[theorem]{Definition}

\newcommand{\bcl}{\begin{center}}
\newcommand{\ecl}{\end{center}}
\newcommand{\brl}{\begin{right}}
\newcommand{\erl}{\end{right}}
\newcommand{\ben}{\begin{enumerate}}
\newcommand{\een}{\end{enumerate}}
\newcommand{\overliner}{\begin{array}}
\newcommand{\earr}{\end{array}}
\newcommand{\btab}{\begin{tabular}}
\newcommand{\etab}{\end{tabular}}
\newcommand{\bdoc}{\begin{document}}
\newcommand{\edoc}{\end{document}}
\newcommand{\beqy}{\begin{eqnarray}}
\newcommand{\eeqy}{\end{eqnarray}}

\newcommand{\beqi}{\begin{eqnarray*}}
\newcommand{\eeqi}{\end{eqnarray*}}
\newcommand{\bitem}{\begin{itemize}}

\newcommand{\eitem}{\end{itemize}}
\newcommand{\nln}{\newline}
\newcommand{\newt}{\newtheorem}

\newcommand{\pa}{\partial}
\newcommand{\re}{{I\!\!R}}
\newcommand{\ren}{\re^N}
\newcommand{\xr}{x\in\re }
\newcommand{\x}{\times}
\newcommand{\dyle}{\displaystyle}
\newcommand{\ene}{{I\!\!N}}
\newcommand{\irn}{\int\limits_{\re^N}}
\newcommand{\io}{\int\limits_{\O}}
\newcommand{\meas}{{\rm meas\,}}

\newcommand{\sign}{{\rm sign}}
\newcommand{\map}{\longrightarrow }
\newcommand{\imp}{\Longrightarrow }
\renewcommand{\div}{\nabla\cdot }
\newcommand{\sen}{{\rm sen\,}}
\newcommand{\tg}{{\rm tg\,}}
\newcommand{\arcsen}{{\rm arcsen\,}}
\newcommand{\arctg}{{\rm arctg\,}}
\newcommand{\supp}{{\textsl supp\ }}
\newcommand{\ity}{\int_{-\iy}^{+\iy}}
\newcommand{\limit}{\lim\limits}
\newcommand{\limi}{\limit_{n\to\infty}}
\newcommand{\sumi}{\sum\limits_{n=1}^{\infty}}
\newcommand{\ulu}{\underline u}
\newcommand{\ulw}{\underline w}
\newcommand{\ulz}{\underline z}
\newcommand{\ulv}{\underline v}
\newcommand{\uls}{\underline s}
\newcommand{\olu}{\overline u}
\newcommand{\olv}{\overline v}
\newcommand{\ols}{\overline s}
\newcommand{\ob}{\overline\b}
\newcommand{\ovar}{\overline\var}
\newcommand{\wv}{\widetilde v}
\newcommand{\wu}{\widetilde u}
\newcommand{\ws}{\widetilde s}
\renewcommand{\a }{\alpha }
\renewcommand{\b }{\beta }
\newcommand{\g }{\gamma}
\newcommand{\G }{\Gamma }
\renewcommand{\d }{\delta }

\newcommand{\D }{\Delta }
\newcommand{\e }{\varepsilon }
\newcommand{\z }{\zeta }
\renewcommand{\l }{\lambda }
\renewcommand{\L }{\Lambda }
\newcommand{\m }{\mu }
\newcommand{\n }{\nabla }
\newcommand{\s }{\sigma }
\newcommand{\Sig }{\Sigma }
\renewcommand{\t }{\tau }
\newcommand{\var }{\varphi }
\renewcommand{\o }{\omega }
\renewcommand{\O }{\Omega }
\newcommand{\bR}{{\bf R}}
\newcommand{\bC}{{\bf C}}
\newcommand{\bZ}{{\bf Z}}
\newcommand{\bN}{{\bf N}}
\newcommand{\bQ}{{\bf Q}}
\newcommand{\bK}{{\bf K}}
\newcommand{\bI}{{\bf I}}
\newcommand{\bv}{{\bf v}}
\newcommand{\bV}{{\bf V}}
\DeclareMathOperator{\suppo}{supp} \DeclareMathOperator{\di}{div}

\def\qed{\unskip\kern 6pt \penalty 500
\raise -2pt\hbox{\vrule \vbox to10pt{\hrule width 4pt
\vfill\hrule}\vrule}\par}

\newenvironment{Proof}{\removelastskip\vskip12pt
plus 1pt \noindent\em\rm}{\hfill {\qed \hskip .2cm}}

\title{Uniqueness in weighted Lebesgue spaces \\
for a class of fractional parabolic and elliptic equations}
\author{Fabio
Punzo\thanks{Dipartimento di Matematica ``F. Enriques",
Universit\`a degli Studi di Milano, via Cesare Saldini 50, 20133
Milano, Italy (fabio.punzo@unimi.it).}\; and Enrico Valdinoci
\thanks{Weierstra{\ss} Institut f\"ur Angewandte Analysis und Stochastik,
Mohrenstra{\ss}e 39, 10117 Berlin, Germany and
Dipartimento di Matematica ``F. Enriques",
Universit\`a degli Studi di Milano,
Via Cesare Saldini 50, 20133 Milano, Italy (enrico@math.utexas.edu).
Supported by the
ERC grant $\epsilon$ ({\it Elliptic Pde's and Symmetry of Interfaces and Layers
for Odd Nonlinearities}). Both authors are supported by PRIN grant 2012 (Critical Point Theory and Perturbative Methods for Nonlinear Differential Equations). }}

\date{}

\begin{document}
\maketitle

\abstract{We investigate uniqueness, in
suitable weighted Lebesgue spaces, of solutions to a class of fractional parabolic and elliptic equations.}


\bigskip
\smallskip

\section{Introduction}\setcounter{equation}{0}
We are concerned with uniqueness of solutions to the following
linear {\it nonlocal} Cauchy problem:
\begin{equation}\label{e1}
\left\{
\begin{array}{ll}
\,   \rho\, \pa_t u + (-\Delta)^{s}u = 0
&\textrm{in}\,\,\re^N\times (0, T]=:S_T
\\& \\
\textrm{ }u \, = 0& \textrm{in\ \ } \re^N\times \{0\} \,,
\end{array}
\right.
\end{equation}
where the coefficient $\rho$, usually referred to as a variable
density, is a positive function only depending on the space
variable $x$ and~$(-\Delta)^s$ denotes the fractional Laplace
operator of order $s\in (0,1)$. Note that when $s=1$ problem
\eqref{e1} has been extensively investigated, and we shall now briefly recall the basics of the classical theory. In particular, it is
well-known that problem
\begin{equation}\label{e2}
\left\{
\begin{array}{ll}
\,   \rho\, \pa_t u -\Delta u = 0 &\textrm{in}\,\,S_T
\\& \\
\textrm{ }u \, = u_0& \textrm{in\ \ } \re^N\times \{0\} \,,
\end{array}
\right.
\end{equation}
where $u_0\in C(\re^N)\cap L^\infty(\re^N),$ admits at most one
bounded solutions if $\rho(x)\to 0$ slowly enough as $|x|\to \infty$
(see, $e.g.$, \cite{EKP}, \cite{IKO}, \cite{KPT} for precise statements). Furthermore,
uniqueness of solutions to problem \eqref{e2}, belonging to
suitable weighted Lebesgue spaces, can be deduced from general
results in \cite{AB}. In fact, for any $\phi\in C(\bar S_T),
\phi>0, p\geq 1$, set
\[L^p_\phi(S_T):= \left\{u: S_T\to \re\,\;\textrm{measurable}\;\,|\, \int_0^T\int_{\re^N}|u|^p\phi(x,t) dx dt<\infty\right\}\,.\]

\smallskip

\noindent {\bf (a) Summary of known results in the local case.}
In \cite{AB}, the operator
\[L u:= \sum_{i,j=1}^N\frac{\pa^2 \big[a_{ij}(x,t) u\big] }{\pa x_i\pa x_j} - \sum_{i=1}^N\frac{\pa \big[b_i(x,t) u\big]}{\pa x_i}+ c(x,t)u - u_t\,\]
is considered; the coefficients of $L$, together with all their
derivatives which appear, are locally bounded functions in $S_T$.
Furthermore, the matrix $A\equiv (a_{ij})$ is assumed to be
positive semidefinite in $\bar S_T$. A function $u\in C(\bar
S_T)$, such that all of its derivatives which appear in $L$ exist
and are locally integrable in $S_T$, is a solution to problem
\begin{equation}\label{e3}
\left\{
\begin{array}{ll}
\,   L\, u = 0 &\textrm{in}\,\,S_T
\\& \\
\textrm{ }u \, = 0& \textrm{in\ \ } \re^N\times \{0\} \,,
\end{array}
\right.
\end{equation}
if equalities in \eqref{e3} are satisfied pointwise. Suppose that
\[ |a_{ij}(x,t)|\leq K_1 ( 1+|x|^2)^{\frac{2-\l}{2}},\; |b_i(x,t)|\leq K_2(|x|^2+1)^{\frac{1}{2}},\;|c(x,t)|\leq K_3(|x|^2+1)^{\frac{\l}2} \]
for almost every $(x,t)\in S_T$, for some constants $\l\geq 0,
K_i>0\; (i=1,2,3)$. In \cite[Theorem 1]{AB} it is shown that if
$u$ is a solution to problem \eqref{e3} and $u\in L^1_{g}(S_T)$,
with
\begin{equation}\label{e63}
g(x)=(|x|^2+1)^{-\a_0} \quad (x\in \re^N)\;\; \textrm{if}\;\;
\l=0,
\end{equation}
or
\begin{equation}\label{e64}
g(x)=e^{-\a_0(|x|^2+1)^{\frac{2-\a}2}}\quad (x\in \re^N)\;\;
\textrm{if}\;\; \l>0\,,
\end{equation}
for some $\a_0>0$, then
\[u\equiv 0\quad \textrm{in}\;\; S_T\,.\]

Let us mention that a crucial step in the proof of this result
consists in the construction of a positive supersolution $\phi\in
C^2(\bar S_T)$ to the adjoint equation
\begin{equation}\label{e4}
\sum_{i,j=1}^N a_{ij}\frac{\pa^2 \phi }{\pa x_i \pa x_j}
+\sum_{i=1}^N b_i\frac{\pa \phi}{\pa x_i} + c \phi + \pa_t\phi = 0
\quad \textrm{in}\;\; S_T.
\end{equation}
Clearly, as a consequence of the previous result one can immediately deduce
uniqueness in $L^1_\psi(S_T)$ for solutions to problem
\[\left\{
\begin{array}{ll}
\,   L\, u = f &\textrm{in}\,\,S_T
\\& \\
\textrm{ }u \, = u_0& \textrm{in\ \ } \re^N\times \{0\} \,,
\end{array}
\right.
\]
where $f$ and $u_0$ are given functions defined in $S_T$ and
$\re^N$, respectively.

\smallskip

Now, suppose that $\rho\in C^2(\re^N).$ If we choose
$$ a_{ij}=\frac 1 {\rho}\d_{ij},\quad  b_i =2 \frac{\pa}{\pa
x_i}\left(\frac 1{\rho}\right),\quad c=\Delta\left(\frac
1{\rho}\right), $$
then, it is easily checked that $u$ is a solution of \eqref{e3} if
and only if it is a solution to \eqref{e1}. Assume that
\begin{equation}\label{e5}
\rho(x)\geq K_1 (1+|x|^2)^{-\frac{\a}2},
\end{equation}
\begin{equation}\label{e6}
\left|\nabla \left(\frac 1{\rho}\right)\right|\le K_2
(1+|x|^2)^{\frac 1 2}, \left|\Delta \left(\frac
1{\rho}\right)\right|\leq K_3(1+|x|^2)^{\frac{2-\a}{2}}
\end{equation} for every $x\in \re^N$, for some $\a\leq 2$. By \cite[Theorem 1]{AB} recalled above, with $\l=2-\a$, $u\equiv 0$ is the unique solution to
problem \eqref{e2} in $L^1_g(S_T)$ with $g$ is defined as in
\eqref{e63}-\eqref{e64}\,.

Moreover, the same conclusion remains true, if we only suppose
$\rho\in C(\re^N)$ instead of $\rho\in C^2(\re^N)$, and we remove
assumption \eqref{e6}. In order to see this, only minor changes in
the proof of \cite[Theorem 1]{AB} are needed. More precisely, the
operator $\hat L u:= \Delta u - \rho \pa_t u$ can play the same
role as $L$; furthermore, we use the fact that, indeed, the
supersolution $\phi$ to equation \eqref{e4}, constructed in
\cite{AB}, is also a supersolution to equation
\[ \rho \pa_t \phi + \Delta \phi =0\quad \textrm{in}\;\; S_T\,.\]

When $\rho\equiv 1$ the results established in \cite{AB} are in
accordance with those in \cite{Ti}, where, not surprisingly, the
same uniqueness class is obtained mainly using the {\it heat
kernel}
\[q(x,y)=\frac 1{(4\pi t)^{\frac N2}} e^{-\frac{|x|^2}{4t}}\quad (x\in \re^N, t>0)\,.\]

\smallskip

\noindent {\bf (b) Summary of known results in the fractional case.}
Observe that also for problem \eqref{e1} some uniqueness results
are known. In fact, in \cite{PT} it is proved that problem
\begin{equation}\label{e1a}
\left\{
\begin{array}{ll}
\,   \rho\, \pa_t u + (-\Delta)^{s}[u^m] = 0 &\textrm{in}\,\, S_T
\\& \\
\textrm{ }u \, = u_0& \textrm{in\ \ } \re^N\times \{0\} \,,
\end{array}
\right.
\end{equation}
$(m\geq 1)$ admits at most one bounded nonnegative solution,
provided condition \eqref{e5} is satisfied, for some $K_1>0,
\a<2s$. Consequently, in particular, $u\equiv 0$ is the unique
solution if $u_0\equiv 0, m=1$. Note that the hypothesis that~$\a<2s$,
when $s=1$ is in agreement with that made above for problem
\eqref{e2}, although for problem \eqref{e2} $\a=2$ was permitted,
too. Moreover, in \cite{BPSV} it is shown that every nonnegative
solution to equation
\[ \pa_t u + (-\Delta)^{s}u=0 \quad \textrm{in}\;\; S_T \]
can be written as
\[u(x,t)=\int_{\re^N} p(x-y, t)u(y,0) dy, \]
where $p$ is the {\it fractional heat kernel} defined by
\[p(x,t):= \frac 1{t^{\frac N{2s}}} P\left(\frac x{t^{\frac 1{2s}}}\right)\quad (x\in \re^N, t>0),\]
\[P(x):= \int_{R^N}e^{ix \cdot \xi-|\xi|^{2s}}d\xi\quad (x\in \re^N)\,. \]
Note that, for some $C>1$, the following two-sided heat kernel
estimate holds (see \cite{BG}):
\begin{equation}\label{e7}
\frac 1 C \min \left\{\frac 1 {t^{\frac N{2s}}}, \frac t{|x|^{N+2s}} \right\} \leq p(x,t)\leq C \min\left\{\frac 1 {t^{\frac N{2s}}}, \frac t{|x|^{N+2s}}\right\} \quad (x\in \re^N, t>0)\,.
\end{equation}
Some uniqueness results for nonlocal parabolic equations have been recently obtained in \cite{AK}, \cite{MP1}, \cite{MP2}.
These works deal in fact with quite general integro-differential operators,
but they require that there exist two constants $C_2>C_1>0, $ such that 
\begin{equation}\label{e200}
C_1\leq \rho \leq C_2\quad \textrm{in}\;\; \re^N\,.
\end{equation}

\smallskip

\noindent{\bf (c) Outline of our results\,.} The main results of this paper will be given in detail
in the forthcoming Theorems \ref{teo1}, \ref{teo3}, \ref{teo4}. We give here a sketchy
outline of these results, describing motivations,
techniques of proofs and differences with the existing literature.

As pointed out for problem \eqref{e2}, we can expect that the
uniqueness class for problem \eqref{e1} with $\rho\equiv 1$ is
related to the fractional heat kernel $p$, and so to its bounds
given in \eqref{e7}. In fact, suppose that condition \eqref{e5} is
satisfied for some $\a\leq 2s$.  We shall prove that the
solution to problem \eqref{e1} is unique in the class
$L^p_{\psi}(S_T)$ with $p\geq 1$ and 
\begin{equation}\label{e66}
\psi(x):=(1+|x|^2)^{-\frac{\b}2}\;\, (x\in \re^N)\,,
\end{equation}
for properly chosen $\b>0$. To be specific, when $\a=0$, we can
choose $\b=N+2s$, in agreement with \eqref{e7}. Furthermore, if
$0<\a\leq 2s$, then $0<\beta<N+2s$, thus the uniqueness class is
smaller; roughly speaking, this is due to the fact that, in this
case, the coefficient $\rho(x)$ makes the diffusion stronger as
$|x|\to\infty$ (see Theorem \ref{teo1} for the precise statement). Observe that, clearly, from such results we can deduce also uniqueness in
$L^p_{\psi}(S_T)$ of solutions to
\begin{equation}\label{e1b}
\left\{
\begin{array}{ll}
\,   \rho\, \pa_t u + (-\Delta)^{s}u = f
&\textrm{in}\,\,S_T
\\& \\
\textrm{ }u \, = u_0& \textrm{in\ \ } \re^N\times \{0\} \,.
\end{array}
\right.
\end{equation}
In order to prove such a uniqueness result we construct a positive supersolution to equation
\begin{equation}\label{e8}
-(-\Delta)^{s} \phi- \rho \pa_t\phi = 0\quad \textrm{in}\;\; S_T\,.
\end{equation}
Indeed, the weight function $\psi$ mentioned above is related to
such a supersolution.

Observe that we also establish similar uniqueness results for the linear elliptic {\it nonlocal} equation
\begin{equation}\label{e9}
(-\Delta)^{s} u + \rho c u \,=\,0 \quad \textrm{in}\;\; \re^N,
\end{equation}
where $c$ is a nonnegative function defined in $\re^N$ (see Theorems \ref{teo3} and \ref{teo4}). Similar
results are stated in \cite{OR} and \cite{P2} for local elliptic
equations in bounded domains, with coefficients that can be
degenerate or singular at the boundary of the domain.

We should observe that the techniques used in the proofs of this paper have several conceptual and technical differences with respect to the classical case of the local equations. On the one hand, we borrow from the classical case the idea of
dealing with the adjoint operator. On the other hand, the classical case
relies on explicit computations based on differentiating some appropriate barriers and cut-off functions and integrating by parts, which are not available in our case. For this reason we have to
perform some ad-hoc integral computations in our case, based on appropriately chosen covering of $\re^N\times \re^N$, some local and global remainder estimates that
rely on some careful paramater adjustments (see Lemmas \ref{lemma1}, \ref{l3}). Furthermore, while in the local case the supersolution to the adjoint equation is defined using the function $g$ defined in \eqref{e63} or \eqref{e64}, in the present situation it is related to the function $\psi$ defined in \eqref{e66} (see the proof of Theorem \ref{teo1}). Moreover, to show that it is indeed a supersolution, clearly, we cannot make explicit computations based on differentiation; instead, we use some properties of the hypergeometric function $_{2}F_1(a,b,c,s)\equiv F(a,b,c,s)$, with
$a,b\in \re, c>0, s\in \re\setminus\{1\}$  (see \cite[Chapters 15.2, 15.4]{DLMF}). Similar computations, for different purposes, when $c>a+b$ have been made in \cite{FerrVerb}; however, we also consider the cases $c=a+b$ and $c>a+b$, that present some differences.  

We should note that whereas the local counterpart of our results established in \cite{AB} only regard the weighted Lebesgue space $L^1_{g}(S_T)$, we can address $L^p_{\psi}(S_T)$, for each $p\geq 1,$ thus such uniqueness results in $L^p_\psi(S_T)$ are new also for $s=1$. Moreover, to the best of our knowledge our results for elliptic equations (see Theorems \ref{teo3} and \ref{teo4}) have not corresponding results in the literature concerning the local case in the whole $\re^N$; some results are only available in bounded domains (see \cite{OR}, \cite{P2}).  Moreover, Theorem \ref{teo3} is proved similarly to the parabolic case, while the proof of Theorem \ref{teo4} is completely new (it does not have a corresponding argument in \cite{OR} or in \cite{P2}). Moreover, it relies on Lemma \ref{l}, which is rather technical. 

Even if, in general, we do not require that the solutions are bounded, as a particular consequence of our uniqueness results it follows that (see Corollary \ref{cor1}) the solution to problem \eqref{e1} is unique in $L^\infty(S_T)$. This generalizes the results in \cite{PT} for linear problems, since in \cite{PT} it was also requested that the solution is bounded and nonnegative; instead, now we do not need any sign condition on the solutions. 
Moreover, there is a substantial difference with uniqueness results in \cite{AK}, \cite{BPSV}, \cite{MP1}, \cite{MP1}. In fact, on the one hand, in \cite{BPSV} only the case $\rho\equiv 1$ is addressed; moreover, the methods used are different form those of the present paper. On the other hand, differently from \cite{AK}, \cite{MP1}, \cite{MP2} we do not make the assumption \eqref{e200}, thus our  density is allowed to vanish
or to be singular as $|x|\to\infty$; moreover, we use completely different techniques. It is worth mentioning that degeneracy or singularity at infinity of the density is very important for the applications, e.g. see for instance, for the local case, \cite{AB}, \cite{EKP}, \cite{IKO}, \cite{OR}\,.
Clearly, the same model with singular or
degenerate density occurs when considering nonlocal diffusion, in
case, for instance, of rarefied media subject to non-Gaussian stochastic processes. 

\medskip

The paper is organized as follows. In Section \ref{mfr} we recall
some preliminaries about fractional Laplacian and we give the
notion of solutions we shall deal with. Then we state our main
results concerning both parabolic and elliptic problems. Section
\ref{pp} is devoted to the proof of results for parabolic
problems, instead those about elliptic equations are proved in
Section \ref{pe}.

\section{Mathematical framework and results}\setcounter{equation}{0} \label{mfr}
The fractional Laplacian $(-\Delta)^{s}$ can be defined by Fourier
transform. Namely, for any function $g$ in the Schwartz class
$\mathcal S$, we say that
$$(-\Delta)^{\sigma/2} g = h\,,$$
if
\begin{equation}
\label{02015}
\hat{h}(\xi)= |\xi|^{\sigma}\hat{g}(\xi).
\end{equation}
Here, we used the notation~$\hat{h}={\mathfrak F} h$ for the Fourier transform of~$h$.
Furthermore, consider the space
\[\mathcal L^s(\re^N):= \left\{u:\re^N\to \re\,\;\textrm{measurable}\, \,|\,\int_{\re^N}\frac{|u(x)|}{1+|x|^{n+2s}}dx<\infty \right\}\,, \]
endowed with the norm
\[ \| u \|_{\mathcal L^s(\re^N)}:=\int_{\re^N}\frac{|u(x)|}{1+|x|^{N+2s}}dx\,.\]
If $u\in \mathcal L^s(\re^N)$ (see \cite{Silv}), then $(-\Delta)^s
u$ can be defined as a distribution, $i.e.$, for any $\varphi\in
\mathcal S$,
\[ \int_{\re^N} \varphi (-\Delta)^s u \, dx \,=\,   \int_{\re^N}  u (-\Delta)^s \varphi \, dx\,. \]
In addition, suppose that, for some
$\g>0$, $u\in \mathcal L^s(\re^N)\cap
C^{2s+\g}(\re^N)$ if $s<\frac 1 2$, or  $u\in \mathcal L^s(\re^N)\cap C^{1,
2s+\g-1}_{loc}(\re^N)\; \textrm{if}\; s\geq \frac 1 2$. Then we have
\begin{equation}
\label{ea1}
(-\Delta)^{s} u(x)=C_{N,s}\,\, \textrm{P.V.}\, \int_{\re^N} \frac{u(x)-u(y)}{|x-y|^{N+2s}}d y\quad (x\in \re^N),
\end{equation}
where $$C_{N,s}=\frac{2^{2s-1}{2s}
\Gamma((N+2s)/2)}{\pi^{N/2}\Gamma(1-s)},$$ $\Gamma$ being the
Gamma function; moreover, $(-\Delta)^s u\in C(\re^N)$. Note that
the constant $C_{N,s}$ satisfies the identity
\[ (-\Delta)^s u = \mathfrak F^{-1} \big( |\xi|^{2s}\mathfrak F u  \big)\,,\quad \xi\in\re^N, u\in \mathcal S\,,\]
so (see \cite{DPVal})
\[C_{N,s}=\left(\int_{\re^N}\frac{1-\cos(\xi_1)}{|\xi|^{N+2s}}d\xi \right)^{-1}\,.\]

\medskip

Concerning the coefficients $\rho$ and $c$, we always make the
following assumptions:
\[
\textrm{\ \ } \left\{
\begin{array}{l}
(i) \quad \;  \rho\in C(\re^N)\,, \rho(x)>0\quad \textrm{for all}\;\; x\in \re^N\,;
\\
(ii)\quad \! \textrm{there exist}\;\, K>0, \a\in \re \;\;\textrm{such that}\;\; \\
\qquad \,\, \rho(x)\geq K \big(1+|x|^2\big)^{-\a}\;\; \textrm{for
all}\;\; x\in \re^N\,;
\end{array}
\right. \leqno(H_0) \]

\[c\in C(\re^N),\; c(x)\geq 0 \quad \textrm{for all}\;\; x\in \re^N\,. \leqno (H_1) \]

Now we can give the definition of solution to problem \eqref{e1} and to equation \eqref{e9}.

\begin{definition}\label{defsoleqp}
We say that a function $u$ is a {\em solution} to equation
\begin{equation}\label{e30}
\rho\, \pa_t u\,+\, (-\Delta)^s u\,=\,0\quad \textrm{in}\;\;
S_T\,,
\end{equation}
if
\begin{itemize}
\item[(i)]  $u\in C(S_T)$, for each $t\in (0,T]$  $u(\cdot, t)\in \mathcal L^s(\re^N)\cap
C^{2s+\g}(\re^N)$ if $s<\frac 1 2$, or  $u(\cdot, t)\in \mathcal L^s(\re^N)\cap C^{1,
2s+\g-1}_{loc}(\re^N)\; \textrm{if}\; s\geq \frac 1 2$,  for some $\g>0$, $\pa_t u\in
C(S_T)$\,;
\item[(ii)] $\rho(x)\pa_t u + C_{N,s} \textrm{P.V}.\, \displaystyle\int_{\re^N} \frac{u(x,t)-u(y,t)}{|x-y|^{N+2s}}d y \,=\,0$\;\;  for all\;\; $(x,t)\in
S_T$\,.
\end{itemize}
Furthermore, we say that $u$ is a {\em supersolution\;
(subsolution)} to equation \eqref{e30}, if in $(ii)$ instead of
$``="$ we have $``\geq"\; (``\leq")$\,.
\end{definition}

\begin{definition}\label{defsolp}
We say that a function $u$ is a solution to problem \eqref{e1} if
\begin{itemize}
\item[(i)] $u\in C(\bar S_T), \pa_t u\in L^1_{loc}(\bar S_T), u\in
L ^1\big((0,T),\mathcal L^s(\re^N)\big)$\,;
\item[(ii)]  $u$ is a solution to equation \eqref{e30}\,;
\item[(iii)] $u(x,0)=0$\;\; for all\;\; $x\in \re^N\,.$
\end{itemize}
\end{definition}

\begin{definition}\label{defsole}
We say that a function $u$ is a solution to equation \eqref{e9} if
\begin{itemize}
\item[(i)]  $u\in \mathcal L^s(\re^N)\cap
C^{2s+\g}(\re^N)$ if $s<\frac 1 2$, or  $u\in \mathcal L^s(\re^N)\cap C^{1,
2s+\g-1}_{loc}(\re^N)\; \textrm{if}\; s\geq \frac 1 2$, for some $\g>0$\,;
\item[(ii)]  $C_{N,s} \textrm{P.V}.\, \displaystyle\int_{\re^N} \frac{u(x)-u(y)}{|x-y|^{N+2s}}d y + \rho(x)c(x)u(x) \,=\,0$\;\;  for all\;\; $x\in \re^N$\,.
\end{itemize}
Furthermore, we say that $u$ is a {\em supersolution\;
(subsolution)} to equation \eqref{e9}, if in $(ii)$ instead of
$``="$ we have $``\geq"\; (``\leq")$\,.
\end{definition}

\subsection{Parabolic equations: results} Next we prove a general criterion for uniqueness of nonnegative solutions to problem
\eqref{e1} in $L^1_{\psi}(S_T),$ where $\psi$ is defined as in \eqref{e66} for some constant $\beta>0$.

\begin{proposition}\label{prop1}
Let assumption $(H_0)-(i)$ be satisfied. Let $u$ be a solution to
problem \eqref{e1} with $|u(\cdot, t)|^p\in \mathcal L^s(\re^N)$ for some
$p\geq 1$, for each $t>0$. Assume that there exists a positive supersolution
$\phi\in C^2(\bar S_T)$ to equation
\begin{equation}\label{e4}
-(-\Delta)^s\phi + \rho \,\pa_t \phi =0 \quad \textrm{in}\;\;
S_T\,,
\end{equation}
such that
\begin{equation}\label{e5}
\phi(x,t)+|\nabla \phi(x,t)|\leq C \psi(x)\quad \textrm{for
all}\;\; (x,t)\in S_T,
\end{equation}
for some constants $C>0$ and $\b>0$. If $u\in L^p_\psi(S_T)$, then
\[u\equiv 0\quad \textrm{in}\;\; S_T\,.\]
\end{proposition}

After having exhibited such a supersolution $\phi$, as a consequence of Proposition \ref{prop1}, we show the following uniqueness theorem.
\begin{theorem}\label{teo1}
Let assumption $(H_0)$ be satisfied. Let $u$ be a solution to
problem \eqref{e1} with $|u(\cdot, t)|^p\in \mathcal L^s(\re^N)$, for some
$p\geq 1$, for each $t>0$. If $u\in L^p_{\psi}(S_T)$, then
\[u\equiv 0 \quad \textrm{in} \;\; S_T ,\]
provided that one of the following conditions holds true:
\begin{itemize}
\item[(i)] $0<\b\leq N-2s, \a\in \re$\,;
\item[(ii)] $N-2s<\b <N, \a\leq 2s$\,;
\item[(iii)] $\b=N, \a<2s$, or instead of $(H_0)-(ii)$, there holds: $\rho(x)\geq K(1+|x|^2)^{-s}\log(1+|x|^2)$ for all $x\in \re^N$, for some $K>0$;
\item[(iv)] $\b>N, \a+\b\leq 2s+N$\,.
\end{itemize}
\end{theorem}

From Theorem \ref{teo1} we deduce the following
\begin{corollary}\label{cor1}
Let assumption $(H_0)$ be satisfied. Let $u$ be a solution to
problem \eqref{e1}. Suppose that $\a<2s$.
\begin{itemize}
\item[(i)] If $$|u(x,t)| \leq C(1+|x|^2)^{\frac{\sigma}2}\quad
\textrm{for all}\;\; x\in S_T,$$ for some $\s\in (0, 2s-\a)$ and
$C>0$, then
$$u\equiv 0\quad \textrm{in}\;\; S_T\,.$$
\end{itemize}
\end{corollary}

In order to prove Corollary \ref{cor1}$(i)$ it suffices to apply
Theorem \ref{teo1} with $\b=N+2s-\a>N$ and $p=1$.

\subsection{Elliptic equations: results} Now we prove a general criterion for uniqueness of nonnegative solutions to
equation \eqref{e9} in $L^1_{\zeta}(\re^N)$. We suppose that there
exists a positive function $\zeta\in C^2(\re^N)$, which solves
\begin{equation}\label{e14}
-(-\Delta)^s\zeta + \rho \,c \,\zeta < 0 \quad \textrm{in}\;\;
\re^N\,.
\end{equation}
Such inequality is meant in the sense that in Definition
\ref{defsole}-$(ii)$, instead of $``="$ we have $``<"$.

\begin{proposition} \label{prop3}
Let assumptions $(H_0)-(i), (H_1)$ be satisfied. Let $u$ be a
solution to equation \eqref{e1b} with $|u|^p\in \mathcal
L^s(\re^N),$ for some $p\geq 1$. Assume that there exists a
positive function $\zeta\in C^2(\re^N)$, which solves \eqref{e6},
and satisfies
\begin{equation}\label{e7}
\zeta(x)+ |\nabla \zeta(x)|\leq C \psi(x)\quad \textrm{for
all}\;\; x\in \re^N,
\end{equation}
for some constants $C>0, \beta>0$\,. If $u\in L^p_\psi(\re^N)$,
then
\[u\equiv 0\quad \textrm{in}\;\; \re^N\,.\]
\end{proposition}
After having exhibited such a supersolution $\zeta$, as a consequence of Propositions \ref{prop3}, we show the following uniqueness theorem.

\begin{theorem}\label{teo3}
Let assumptions $(H_0)-(H_1)$ be satisfied. Let $u$ be a solution
to equation \eqref{e1b} with $|u|^p\in \mathcal L^s(\re^N),$ for
some $p\geq 1$. Suppose that, for some $c_0>0,$
\begin{equation}\label{e8}
c(x)\geq c_0\quad \textrm{for all}\;\;x\in \re^N\,.
\end{equation}
If $u\in L^p_{\psi}(\re^N)$, then
\[u\equiv 0 \quad \textrm{in} \;\; \re^N ,\]
provided $\a, \b$ satisfy the same conditions as in Theorem
\ref{teo1}, and $p c_0K$ is large enough when $(ii)$ or $(iii)$ or
$(iv)$ holds true.
\end{theorem}

Analogously to Corollary \ref{cor1}, we have the following
\begin{corollary}\label{cor2}
Let assumptions $(H_0)-(H_1)$ be satisfied. Let $u$ be a solution
to equation \eqref{e1b} with $|u|^p\in \mathcal L^s(\re^N),$ for
some $p\geq 1$ . Suppose that $pc_0K$ is large enough and $\a<2s$.
\begin{itemize}
\item[(i)] If $$|u(x)| \leq C(1+|x|^2)^{\frac{\sigma}2}\quad
\textrm{for all}\;\; x\in \re^N,$$ for some $\s\in (0, 2s-\a)$ and
$C>0$, then
$$u\equiv 0\quad \textrm{in}\;\; \re^N\,.$$
\end{itemize}
\end{corollary}

\begin{remark}
{\rm The hypothesis $p c_0 K$ large enough made in Theorem
\ref{teo3} and in Corollary \ref{cor2} will be specified in the
proof of Theorem \ref{teo3}.}
\end{remark}

\medskip

Let us now introduce the Riesz kernel of the $s$-Laplacian:
$$ I_{2s}:=\frac {k_{N,s}}{|x|^{N-2s}} \ \ \ \forall x \in \re^N \, , $$
where $k_{N,s}$ is a suitable positive constant only depending on
$s$ and $N$.

Let $F\in  C^\infty_c(\re^N)$ with $F \geq 0, F\not\equiv 0$.
Define
\begin{equation}\label{e76c}
\phi = I_{2s} \ast F \quad \textrm{in}\;\; \re^N\,.
\end{equation}
Clearly, we have that
\begin{equation}\label{e6c}
(-\Delta)^s \phi =  F \quad \textrm{in} \ \re^N \,.
\end{equation}
Furthermore, it is easily checked that $ \phi \in C^\infty(\re^N)$
and, for some $0<C_0< C_1$,
\begin{equation}\label{e7c}
\frac{C_0}{1+|x|^{N-2s}}\leq \phi(x) + \left|\nabla \phi(x)
\right| \leq\frac{C_1}{1+ |x|^{N-2s}} \quad \textrm{for all}\;\; x
\in \re^N\,.
\end{equation}

If $u\in L^p_{c\rho\phi}(\re^N)$, then we can drop the request $p
c_0 K $ big enough made in Theorem \ref{teo3}. This is the content
of the next result, which will be proved by different methods from
those used to prove Theorem \ref{teo3}.

\begin{theorem}\label{teo4}
Let assumptions $(H_0)-(H_1)$ be satisfied. Let $u$ be a solution
to equation \eqref{e1b}. Suppose that condition \eqref{e8} is
satisfied for some $c_0>0.$ If $u\in L^p_{c\rho\phi}(\re^N)$ for
some $p\geq 1$, then
\[u\equiv 0 \quad \textrm{in} \;\; \re^N\,.\]
\end{theorem}

As a consequence of Theorem \ref{teo4} and \eqref{e7c} we
immediately get the next result.
\begin{corollary}\label{cor4}
Let assumptions $(H_0)-(H_1)$ be satisfied. Let $u$ be a solution
to equation \eqref{e1b}. Suppose that $c\in L^\infty(\re^N)$ and
that condition \eqref{e8} is satisfied for some $c_0>0$. If $u\in
L^p_{1+|x|^{N-2s+\a}}(\re^N)$, then
\[u\equiv 0 \quad \textrm{in} \;\; \re^N\,.\]
\end{corollary}

\section{Parabolic equations: proofs}\setcounter{equation}{0} \label{pp}
\subsection{Preliminary results}
This Subsection is devoted to some preliminary results that will be used in the sequel.
To begin with, let us observe that if $f,g  \in \mathcal L^s(\re^N)\cap
C^{2s+\g}(\re^N)$ if $s<\frac 1 2$, or  $f,g\in \mathcal L^s(\re^N)\cap C^{1,
2s+\g-1}_{loc}(\re^N)\; \textrm{if}\; s\geq \frac 1 2$, for some $\g>0$, and $fg\in \mathcal
L^s(\re^N)$, then it is easily checked that
\begin{equation}\label{e20}
(-\Delta)^s[f(x)g(x)]= f(x)(-\Delta)^{s}g(x) + g(x)(-\Delta)^sf(x)
- \mathcal B(f,g)(x)\quad \textrm{for all}\;\; x\in \re^N\,,
\end{equation}
where $\mathcal B(f,g)$ is the bilinear form given by
\[\mathcal B(f,g)(x):=C_{N,s}\int_{\re^N}\frac{[f(x)-f(y)][g(x)-g(y)]}{|x-y|^{N+2s}}dy \quad \textrm{for all}\;\; x\in \re^N\,.\]
Take a cut-off function $\g\in C^\infty([0,\infty)), 0\leq \g\leq
1$ with
\begin{equation}\label{e81}
\g(r)=\left\{
\begin{array}{ll}
\,  1 &\textrm{if}\,\,0\leq r\leq \frac 1 2
\\& \\
\textrm{ } 0& \textrm{if\ \ } r\geq 1 \,;
\end{array}
\right.
\end{equation}
for any $R>0$ let
\begin{equation}\label{e82}
\g_R(x):= \g\left(\frac{|x|}{R}\right) \quad \textrm{for all}\;\; x\in \re^N\,.
\end{equation}
For any $\t\in (0,T)$ let
\[S_\t:= \re^N\times (0,\t]\,.\]
We shall use next
\begin{lemma}\label{lemma1}
Let $\t\in (0,T), \phi\in C^2(\bar S_\t), \phi>0$; suppose that \eqref{e5} is satisfied. Let $u\in
L^1_{\psi}(S_\t)$. Then
\[\int_0^\t\int_{\re^N}|u(x,t)|\phi(x,t)|(-\Delta)^s\g_R(x)|dx dt +\int_0^\t\int_{\re^N}|u(x,t)|\,|\mathcal B(\phi, \g_R)(x)|dx dt\to 0 \]
as $R\to \infty\,.$
\end{lemma}
Observe that a similar result was obtained in the proof Theorem 2.1 in \cite{BPSV}. However, in \cite{BPSV} different hypotheses were made. To be specific,
it was assumed that $u\in L^1\big((0,T); \mathcal L^s(\re^N)\big)$, and
moreover that, for some $C>0$,
$$\phi(x,t) + |\nabla \phi(x,t)|\leq \frac{C}{1+|x|^{N+2s}}\quad \textrm{for all}\;\; (x,t)\in S_T.$$
Observe that our proof require various quite important changes. In particular, we need to use a convenient covering of $\re^N\times \re^N$ that is a little different from that in \cite{BPSV}; moreover, we shall use different estimates in some regions of $\re^N\times \re^N$. 

\medskip

\noindent{\it Proof\,.\,\,} Observe that, for all $x\in \re^N$,
\begin{equation}\label{e80}
|(-\Delta)^{s}\gamma_{R}(x)|=R^{-2s }\left|
\left((-\Delta)^{s}\gamma\right)\left(\frac{x}{R}\right)\right|\leq C R^{-2s}.
\end{equation}
Since $u\in L^1_\psi(S_T)$ and \eqref{e5} holds, from \eqref{e80} it follows that
$$
 \int_0^\t\int_{\re^N}|u(x,t)| \phi(x,t)|(-\Delta)^s\g_R(x)|dx dt\leq C R^{-2s }\int_{0}^{\t}{\int_{\re^N}{|u(x,t)|\phi(x,t)\, dx\, dt}}\leq C R^{- 2s},
$$
so
\begin{equation}\label{step3}
\lim_{R\to\infty}{  \int_0^\t\int_{\re^N}|u(x,t)|\phi(x,t)|(-\Delta)^s\g_R(x)|dx dt    }=0.
\end{equation}

\smallskip

Now we are going to estimate
\begin{equation}\label{e89}
I(R):= \int_0^\t\int_{\re^N}|u(x,t)|\,|\mathcal B(\phi, \g_R)(x)|dx dt\,
\end{equation}
To do this, we cover $\re^N\times\re^N$ with six domains. In fact,
$$\re^{2N}=\left(\bigcup_{k=1}^{5}{A_{k}}\right)\cup \mathcal C,$$
where
$$A_1:=\{(x,y):\, |x|> R/2,\, |y|\leq R/8\},\quad A_2:=\{(x,y):\, |x|\leq R/8,\, |y|> R/2\},$$
$$A_3:=\{(x,y):\, |x|\geq 2R,\, R/8<|y|< R\},\quad A_4:=\{(x,y):\, R/8<|x|< R,\, |y|\geq 2R \},$$
$$A_5:=\{(x,y):\, R/8<|x|< 2R,\, R/8<|y|< 2R\}$$
and
$$\mathcal C:=\{(x,y):\, |x|\leq R/2, \,|y|\leq R/2\}\cup \{(x,y):\, |x|\geq R, |y|\geq R\}.$$

\smallskip

This covering of $\re^N\times \re^N$ is represented (for $N=1$) in the following picture:
\begin{center}
\begin{tikzpicture}
\def\R{2}
\coordinate (O) at (0,0);

\fill [color=black!20] (O)--(\R/2,0)--(\R/2,\R/2)--(0,\R/2)--cycle;

\fill [color=black!20] (\R,3*\R)--(\R,\R)--(3*\R,\R)--(3*\R,3*\R)--cycle;

\draw (-1,0)--(3*\R,0);
\draw (0,-1)--(0,3*\R);
\draw (O)--(3*\R,3*\R);
\draw (0,\R/2)--(\R/2,\R/2)--(\R/2,0);
\draw (\R/2,\R/8)--(3*\R,\R/8);
\draw (2*\R,\R/8)--(2*\R,\R);
\draw (\R,3*\R)--(\R,\R)--(3*\R,\R);
\draw (\R/8,\R/2)--(\R/8,3*\R);
\draw (\R/8,2*\R)--(\R,2*\R);

\draw ($(\R/8,0)+(0,-\R/25)$)--($(\R/8,0)+(0,+\R/25)$);
\draw (\R/8,0) node [below,yshift=-4pt] {$\frac R8$};

\draw ($(\R/2,0)+(0,-\R/25)$)--($(\R/2,0)+(0,+\R/25)$);
\draw (\R/2,0) node [below,yshift=-4pt] {$\frac R2$};

\draw ($(\R,0)+(0,-\R/25)$)--($(\R,0)+(0,+\R/25)$);
\draw (\R,0) node [below,yshift=-4pt] {$R$};

\draw ($(2*\R,0)+(0,-\R/25)$)--($(2*\R,0)+(0,+\R/25)$);
\draw (2*\R,0) node [below,yshift=-4pt] {$2R$};

\draw ($(0,\R/8)+(-\R/25,0)$)--($(0,\R/8)+(\R/25,0)$);
\draw (0,\R/8) node [left,xshift=-4pt] {$R/8$};

\draw ($(0,\R/2)+(-\R/25,0)$)--($(0,\R/2)+(\R/25,0)$);
\draw (0,\R/2) node [left,xshift=-4pt] {$R/2$};

\draw ($(0,\R)+(-\R/25,0)$)--($(0,\R)+(\R/25,0)$);
\draw (0,\R) node [left,xshift=-4pt] {$R$};

\draw ($(0,2*\R)+(-\R/25,0)$)--($(0,2*\R)+(\R/25,0)$);
\draw (0,2*\R) node [left,xshift=-4pt] {$2R$};

\draw (3*\R/2,\R/16) node {$A_1$} (\R/16,3*\R/2) node {$A_2$} (5*\R/2,\R/2) node {$A_3$} (\R/2,5*\R/2)
node {$A_4$} (\R/2,\R) node {$A_5$} (\R,\R/2) node {$A_5$}  (3*\R/2,5*\R/2) node {$\mathcal C$}
(5*\R/2,3*\R/2) node {$\mathcal C$}  (3*\R/8,\R/8) node {$\mathcal C$}   (\R/8,3*\R/8) node {$\mathcal C$};

\end{tikzpicture}
\end{center}

{F}rom \eqref{e81} and~\eqref{e82}, we have that $\g_R(x)-\g_R(y)=0$ if $(x,y)\in \mathcal C$,
and so
\begin{eqnarray}
I(R)&=&\int_{0}^{t_0}{\int_{\re^N}{|u(x,t)|
\int_{\re^N}{\frac{|\phi(x,t)-\phi(y,t)|\,|
\g_{R}(x)-\g_{R}(y)|}{|x-y|^{N+2s}}\, dy}\, dx\, dt}}\nonumber\\
&\le&\sum_{k=1}^{5}{I^{A_k}(R)},\label{step4}
\end{eqnarray}
where
$$I^{A_k}(R)=\int_{0}^{t_0}{\int_{A_k}{|u(x,t)|\frac{|
\phi(x,t)-\phi(y,t)|\,|\g_{R}(x)-\g_{R}(y)|}{|x-y|^{N+2s}}\, dy}\, dx\, dt},$$
for $k=1,\ldots,5$.

We are going to estimate each of these five integral separately. For all $(x,y)\in A_1$ we get
\begin{equation}\label{e83}
|x-y|\geq C |x|\,,
\end{equation}
and
\begin{equation}\label{e84}
|x-y|\geq\frac{R}{2}-\frac{R}{8} \geq  \frac R 4+|y|\,.
\end{equation}
Hence
\begin{equation}\label{e86}
\frac 1{|x-y|^{N+2s}}\leq \frac C{|x|^{\b}\big( \frac R 4+ |y|\big)^{N+2s-\b}}\,.
\end{equation}
Moreover, from \eqref{e5} it follows that, if~$(x,y)\in A_1$, then
\begin{equation}\label{e85}
|\phi(x,t)|+|\phi(y,t)|\leq \frac{
C}{1+|y|^{\b}}.
\end{equation}
Inequalities \eqref{e86} and \eqref{e85} yield, for any $R>1$,
\begin{eqnarray}
I^{A_1}(R)
&\leq&
 C \int_{0}^{\t}{\int_{|x|> R/2}{\frac{|u(x,t)|}{|x|^{\b}} \int_{|y|\leq R/8}{\frac{1}{1+|y|^{N+2s}}\, dy}\, dx\, dt}}\nonumber\\
&\leq&C\int_{0}^{\t}{\int_{|x|> R/2}{\frac{|u(x,t)|}{|x|^{\b}}\, dx\, dt}}.\label{step5}
\end{eqnarray}

\smallskip

For all $(x,y)\in A_2$ we have
\begin{equation}\label{nuevo}
|\phi(x,t)|+|\phi(y,t)|\leq \frac{
C}{1+|x|^{\b}},
\end{equation}
and
\begin{equation}\label{e84b}
|x-y| \geq C|y|\,.
\end{equation}
In view of \eqref{nuevo} and \eqref{e84b}, we obtain
\begin{eqnarray}
I_{2}^{A_2}(R)
&\leq&
\int_{0}^{\t}{\int_{|x|\leq R/8}{\frac{|u(x,t)|}{1+ |x|^{\b}}\int_{|y|>R/2}{\frac{C}{|y|^{N+2s}}\, dy}\, dx\, dt}}\nonumber\\
&\leq&CR^{-2s}\int_{0}^{\t}{\int_{|x|\le R/8}{\frac{|u(x,t)|}{1+|x|^{\b}}\, dx\, dt}}.\label{step6}
\end{eqnarray}

\smallskip
Also, for all $(x,y)\in A_3$ we have that \eqref{e83}, \eqref{e85} and \eqref{e84b} hold true. From \eqref{e83} and \eqref{e84b} we get
\begin{equation}\label{e87}
\frac 1{|x-y|^{N+2s}}\leq \frac C{|x|^\b|y|^{N+2s-\b}}\,.
\end{equation}
So, due to \eqref{e85} and \eqref{e87}, we obtain
\begin{eqnarray}\label{step7}
I^{A_3}(R)\leq C \int_0^\t \int_{|x|\geq 2R}\frac{|u(x,t)|}{|x|^\b} \int_{\frac R 8<|y|<R}\frac 1{|y|^{N+2s}} dy\leq \nonumber\\
\leq CR^{-2s}\int_{0}^{\t}{\int_{|x|\geq 2R}{\frac{|u(x,t)|}{|x|^{\b}}\, dx\, dt}}.
\end{eqnarray}
For all
$(x,y)\in A_4$, we have that  \eqref{nuevo} and \eqref{e84b} hold true. Hence,
\begin{eqnarray}
I^{A_4}(R)
&\leq&
\int_{0}^{\t}{\int_{R/8\leq |x|\leq R}{\frac{|u(x,t)|}{1+ |x|^{\b}}\int_{|y|>2 R}{\frac{C}{|y|^{N+2s}}\, dy}\, dx\, dt}}\nonumber\\
&\leq&CR^{-2s}\int_{0}^{\t}{\int_{R/8\leq |x|\le R}{\frac{|u(x,t)|}{1+|x|^{\b}}\, dx\, dt}}.\label{step8}
\end{eqnarray}
Then, using the Monotone Convergence Theorem, since $u\in L^1_\psi(S_T)$, from \eqref{step5}, \eqref{step6}, \eqref{step7} and \eqref{step8} it follows that
\begin{equation}\label{IA1234}
\lim_{R\to\infty}I^{A_1}(R)=\lim_{R\to\infty}I^{A_2}(R)=\lim_{R\to\infty}
I^{A_3}(R)=\lim_{R\to\infty}I^{A_4}(R)=0.
\end{equation}
To estimate $I^{A_5}(R)$ we will consider separately the cases $s\in\left(0,\frac 12\right)$ and $s\in\left [\frac 12 ,1\right)$.

\

Let $s\in\left(0,\frac 12\right)$ and $(x,y)\in A_5$. Since in $A_5$
the roles of $x$ and $y$ are symmetric, from \eqref{e5} we can infer that
\begin{equation}\label{b}
|\phi(x,t)|+|\phi(y,t)|\leq \frac{C}{|x|^{\b}}.
\end{equation}
Furthermore,
\begin{equation}\label{e88}
|\g_{R}(x)-\g_{R}(y)|\leq \frac{C}{R}|x-y|.
\end{equation}
Hence
\begin{equation}\label{step9}
I^{A_5}(R)\leq \frac{C}{R}\int_{0}^{\t}{\int_{\frac{R}{8}\leq|x|\leq 2R}{\frac{|u(x,t)|}{|x|^{\b}}\int_{\frac{R}{8}\leq|y|\leq 2R}{\frac{1}{|x-y|^{N+2s-1}}\, dy}\, dx\, dt}}.
\end{equation}
By the change of variables $\tilde y:=x-y$, since $s\in \left(0, \frac 12\right)$, from \eqref{step9} it follows  that
\begin{eqnarray}
I^{A_5}(R)
&\leq&
\frac{C}{R}\int_{0}^{\t}{\int_{\frac{R}{8}\leq|x|\leq 2R}{\frac{|u(x,t)|}{|x|^{\b}}\int_{\frac{R}{8}\leq|x-\tilde y|\leq 2R}{\frac{1}{|\tilde y|^{N+2s-1}}\, d\tilde y}\, dx\, dt}}\nonumber\\
&\leq&\frac{C}{R}\int_{0}^{\t}{\int_{\frac{R}{8}\leq|x|\leq 2R}{\frac{|u(x,t)|}{|x|^{\b}}\int_{|\tilde y|\leq 4R}{\frac{1}{|\tilde y|^{N+2s-1}}\, d\tilde y}\, dx\, dt}}.\nonumber\\
&\leq&
CR^{-2s}\int_{0}^{\t}{\int_{\frac{R}{8}\leq|x|\leq 2R}{\frac{|u(x,t)|}{|x|^{\b}}\, dx\, dt}}.\label{step10}
\end{eqnarray}
Therefore, since $u\in L^1_\psi(S_T)$,  we conclude that
\begin{equation}\label{step11}
\lim_{R\to\infty}I^{A_5}(R)=0,\, \qquad{\mbox{ when }}s\in\left (0,\frac 12\right).
\end{equation}
Now, let $s\in\left[\frac 12 ,1\right)$. By \eqref{e5}, we get
\begin{equation}\label{step111}
|\phi(x,t)-\phi(y,t)|\leq\frac{C}{1+|z|^{\b}}|x-y|,
\end{equation}
for some $z$  in the segment joining $x$ and $y$.
For any $R>0$ let
$$Q_R\equiv Q:=\left\{(x,y)\in A_5:\, |x-y|\leq \frac{R}{100}\right\}.$$
Note that, if $(x,y)\in Q$ then every point $z$ lying on the segment from $x$
to $y$ verifies $|z|\geq C|x|$. Hence, since $s\in \left[ \frac 1 2,1\right)$,  \eqref{step111} and \eqref{e88} yield
\begin{eqnarray}
&&\int_{0}^{\t}{\int_{(x,y)\in Q}{|u(x,t)|\frac{|(\g_{R}(x)-\g_{R}(y))(\phi(x,t)-\phi(y,t))|}{|x-y|^{\b}}\, dy\, dx\, dt}}\nonumber\\
&\leq&\int_{0}^{\t}{\int_{(x,y)\in Q}{|u(x,t)|\frac{C}{R|x|^{\b}|x-y|^{N+2s-2}}\, dy\, dx\, dt}}\nonumber\\
&\leq&\frac{C}{R}\int_{0}^{\t}{\int_{\frac{R}{8}\leq|x|\leq 2R}{\frac{|u(x,t)|}{|x|^{\b}}\int_{\frac{R}{8}\leq|y|\leq 2R}{\frac{1}{|x-y|^{N+2s-2}}\, dy}\, dx\, dt}}\nonumber\\
&\leq&\frac{C}{R}\int_{0}^{\t}{\int_{\frac{R}{8}\leq|x|\leq 2R}{\frac{|u(x,t)|}{|x|^{\b}}\int_{|\tilde y|\leq 4 R}{\frac{1}{|\tilde y|^{N+2s-2}}\, d\tilde y}\, dx\, dt}}\nonumber\\
&\leq&CR^{1-2s}\int_{0}^{\t}{\int_{\frac{R}{8}\leq|x|\leq 2R}{\frac{|u(x,t)|}{|x|^{\b}}\, dx\, dt}}.\label{step12}
\end{eqnarray}
On the other hand, if 
$(x,y)\in A_5\setminus Q$ we have that
\begin{equation}\label{c}
|x-y|>\frac{R}{100}\geq C|y|.
\end{equation}
Then, by \eqref{b} and \eqref{c},
\begin{eqnarray}
&&\int_{0}^{\t}{\int_{(x,y)\in A_5\setminus Q}{|u(x,t)|\frac{|(\g_{R}(x)-\g_{R}(y))(\phi(x,t)-\phi(y,t))|}{|x-y|^{N+2s}}\, dy\, dx\, dt}}\nonumber\\
&\leq&\frac{C}{R}\int_{0}^{\t}{\int_{(x,y)\in A_5\setminus Q}{|u(x,t)|\frac{|(\phi(x,t)-\phi(y,t))|}{|x-y|^{N+2s-1}}\, dy\, dx\, dt}}\nonumber\\
&\leq&\frac{C}{R}\int_{0}^{\t}{\int_{(x,y)\in A_5\setminus Q}{\frac{|u(x,t)|}{|x|^{\b}}\frac{1}{|y|^{N+2s-1}}\, dy\, dx\, dt}}\nonumber\\
&\leq&\frac{C}{R}\int_{0}^{\t}{\int_{\frac{R}{8}\leq|x|\leq 2R}{\frac{|u(x,t)|}{|x|^{\b}}\int_{\frac{R}{8}\leq|y|\leq 2R}{\frac{1}{|y|^{N+2s-1}}\, dy}\, dx\, dt}}\nonumber\\
&\leq&CR^{-2s}\int_{0}^{\t}{\int_{\frac{R}{8}\leq|x|\leq 2R}{\frac{|u(x,t)|}{|x|^{\b}}\, dx\, dt}}.\label{step13}
\end{eqnarray}
Therefore, from \eqref{step12} and \eqref{step13} we have
\begin{eqnarray}
I^{A_5}(R)&\leq&\int_{0}^{\t}{\int_{(x,y)\in Q}{|u(x,t)|\frac{|(\g_{R}(x)-\g_{R}(y))(\phi(x,t)-\phi(y,t))|}{|x-y|^{N+2s}}\, dy\, dx\, dt}}\nonumber\\
&+&\int_{0}^{\t}{\int_{(x,y)\in A_5\setminus Q}{|u(x,t)|\frac{|(\g_{R}(x)-\g_{R}(y))(\phi(x,t)-\phi(y,t))|}{|x-y|^{N+2s}}\, dy\, dx\, dt}}\nonumber\\
&\leq&CR^{1-2s}\int_{0}^{\t}{\int_{\frac{R}{8}\leq|x|\leq 2R}{\frac{|u(x,t)|}{|x|^{\b}}\, dx\, dt}}\nonumber\\
&+&CR^{-2s}\int_{0}^{\t}{\int_{\frac{R}{8}\leq|x|\leq 2R}{\frac{|u(x,t)|}{|x|^{\b}}\, dx\, dt}}.
\end{eqnarray}
Then, since $u\in L^1_\psi(S_T)$, using the Monotone Convergence Theorem we obtain
\begin{equation}\label{step14}
\lim_{R\to\infty} I^{A_5}(R)=0,\, \qquad{\mbox{ when }}s\in\left[\frac 12 ,1\right).
\end{equation}
That is, by \eqref{step11} and \eqref{step14}, we get
\begin{equation}\label{IA5}
\lim_{R\to\infty} I^{A_5}(R)=0,\, \qquad{\mbox{ whenever }} s\in(0,1).
\end{equation}
Putting together \eqref{step4},  \eqref{IA1234} and \eqref{IA5} it follows that
\begin{equation}\label{step15}
\lim_{R\to\infty}I(R)=0,\, \qquad{\mbox{ when }}s\in(0,1).
\end{equation}
From \eqref{step3}, \eqref{e89} and \eqref{step15} the conclusion follows. \hfill $\square$

\bigskip

\subsection{Proof of Proposition \ref{prop1}}
The next lemma will be used.
\begin{lemma}\label{l3}
Let $G\in C^2(\re; \re)$ be a convex function. Let $u\in \mathcal
L^s(\re^N)\cap C^{2s+\g}(\re^N)$ if $s<\frac 1 2$, or $u\in
\mathcal L^s(\re^N)\cap C^{1, 2s+\g-1}_{loc}(\re^N)\;
\textrm{if}\; s\geq \frac 1 2$,  for some $\g>0$. Suppose that
$G(u)\in \mathcal L^s(\re^N)$. Then
\begin{equation}\label{e12}
(-\Delta)^s[G(u)]\leq G'(u)(-\Delta)^s u\quad
\textrm{in}\,\;\re^N\,.
\end{equation}
\end{lemma}

\noindent{\it Proof\,.} We can choose, by a suitable convolution, a sequence $\{u_n\}\subset
\mathcal S$ uniformly bounded in $C^{2s+\g}_{loc}(\re^N)$, if
$s<\frac 1 2$, or in $C^{1, 2s+\g-1}_{loc}(\re^N)\; \textrm{if}\;
s\geq \frac 1 2$, for some $\g>0$, with $u_n\to u$ as $n\to\infty$
both in $\mathcal L^s(\re^N)$ and locally uniformly in $\re^N$.
Since $G\in C^2(\re;\re)$ and $G(u)\in \mathcal L^s(\re^N)$,
analogously to the proof of \cite[Proposition 2.1.4]{Silv} we have
that
\[(-\Delta)^s u_n\to (-\Delta)^s u, (-\Delta)^s[G(u_n)]\to (-\Delta)^s G(u)\quad \textrm{as}\,\, n\to\infty, \]
locally uniformly in $\re^N$. From \cite[Lemma 4.1]{FFV} we have
\[(-\Delta)^s[G(u_n)]\leq G'(u_n)(-\Delta)^s u_n\quad
\textrm{in}\,\;\re^N\,.\] So, passing to the limit as $n\to
\infty$ we get \eqref{e12}.

\smallskip
\medskip

\noindent{\it Proof of Proposition \ref{prop1}}\,. Let $\t\in
(0,T)$. Take a nonnegative function $v\in C^2(\bar S_\t)$ with
$supp\;\, v(\cdot, t)$ compact for each $t\in [0,\t]$. Moreover,
take a function $w\in C(\bar S_\t)$ such that for each $t\in
(0,\t]$, $w(\cdot, t)\in \mathcal L^s(\re^N)\cap C^{2s+\g}(\re^N)$
if $s<\frac 1 2$, or  $w(\cdot, t)\in \mathcal L^s(\re^N)\cap
C^{1, 2s+\g-1}_{loc}(\re^N)\; \textrm{if}\; s\geq \frac 1 2$, for
some $\g>0, w\in L^1\big((0,\t);\mathcal L^s(\re^N)\big)$. For any
$\epsilon\in (0,\tau)$, integrating by parts we have:
\begin{equation}\label{e13}
\begin{split}
 \int_0^\t\int_{\re^N} v \big[-(-\Delta)^s w- \rho \pa_t w \big]\,  dx dt =
 \int_0^\t \int_{\re^N} w\big[- (- \Delta)^s v+\rho\pa_t v\big]\, dx dt\\
  -\int_{\re^N}\rho(x) v(x,\t)w(x,\t) dx +\int_{\re^N}\rho(x) v(x,\epsilon) w(x,\epsilon) dx\,.
\end{split}
\end{equation}
Let $p\geq 1.$ For any $\a>0$, set
\begin{equation}\label{e23}
G_\a(r):=(r^2+\a)^{\frac p 2} \quad \textrm{for all}\;\; r\in
\re\,.
\end{equation}
It is easily seen that
\begin{equation}\label{e14}
G_\a''(r)\geq 0 \quad \textrm{for all}\;\; r\in\re\,.
\end{equation}
By \eqref{e1},
\begin{equation}\label{e15}
\rho\pa_t[G_\a(u)]=G_\a'(u)\pa_t u=-G_\a'(u)(-\Delta)^s u\quad
\textrm{for all}\;\; (x,t)\in S_T\,.
\end{equation}
From \eqref{e14}, \eqref{e15} and Lemma \ref{l3} we obtain
\begin{equation}\label{e16}
\rho \pa_t[G_\a(u)] + (-\Delta)^s[G_\a(u)]\leq 0 \quad
\textrm{in}\;\; S_T\,.
\end{equation}
So, from \eqref{e13} with $w=G_\a(u)$ and \eqref{e15} we obtain
\begin{equation}\label{e17}
\begin{split}
\int_{\re^N} \rho(x) G_\a[u(x,\t)] v(x,\t) dx \leq  \int_0^\t
\int_{\re^N} G_\a(u) \big[- (- \Delta)^s v+\rho\pa_t v\big] dx dt
\\ + \int_{\re^N}\rho(x) v(x,\epsilon) G_\a[u(x,\epsilon)] dx\,.
\end{split}
\end{equation}
Letting $\epsilon\to 0^+$ in \eqref{e17}, by the dominated
convergence theorem,
\begin{equation}\label{e18}
\begin{split}
\int_{\re^N} \rho(x) G_\a[u(x,\t)] v(x,\t) dx \leq  \int_0^\t
\int_{\re^N} G_\a(u) \big[- (- \Delta)^s v+\rho\pa_t v\big] dx dt
\\ + \a \int_{\re^N}\rho(x) v(x,0) dx\,.
\end{split}
\end{equation}
Now, letting $\a\to 0^+$, by the dominated convergence theorem,
\begin{equation}\label{e19}
\int_{\re^N} \rho(x) |u(x,\t)|^p v(x,\t) dx \leq  \int_0^\t
\int_{\re^N} |u|^p \big[- (- \Delta)^s v+\rho\pa_t v\big] dx dt
\end{equation}

For any $R>0$, we can choose
\[v(x,t):= \phi(x,t)\g_R(x)\quad \textrm{for all}\;\; (x, t)\in \bar S_\t\,.\]
From \eqref{e19}, using the fact that $\phi$ is a supersolution to
equation \eqref{e4} and $\g_R\geq 0$, we obtain
\begin{equation}\label{e20}
\begin{split}
-(-\Delta)^s v+\rho \pa_t v=\gamma_R \big[-(-\Delta)^s \phi + \rho
\pa_t \phi\big]- \phi(-\Delta)^s\g_R + \mathcal B(\phi,\g_R)
\\\leq - \phi(-\Delta)^s\g_R + \mathcal B(\phi,\g_R) \quad
\textrm{in}\;\; S_\t\,.
\end{split}
\end{equation}
Since $|u|^p\geq 0$, by \eqref{e19} and \eqref{e20} we conclude
that
\begin{equation}\label{e21}
\begin{split}
\int_{\re^N} \rho(x) |u(x,\t)|^p \phi(x,\t)\g_R(x) dx \leq
\int_0^\t \int_{\re^N}
|u|^p \big[ - \phi(-\Delta)^s\g_R + \mathcal B(\phi,\g_R) \big]dx dt\\
\leq \int_0^\t \int_{\re^N}|u|^p\big[ |\phi||(-\Delta)^s\g_R| +
|\mathcal B(\phi,\g_R)|\big]dx dt\,.
\end{split}
\end{equation}
Hence, from Lemma \ref{lemma1} with $v=|u|^p$ and the monotone
convergence theorem, sending $R\to \infty$ in \eqref{e21} we get
\begin{equation}\label{e22}
\int_{\re^N}\rho(x) |u(x,\t)|^p\phi(x,\t)dx\leq 0\,.
\end{equation}
From \eqref{e22}, $(H_0)-(i)$, since $\phi>0$ in $S_\t$ and
$|u|^p\geq 0$ we infer that $u\equiv 0$ in $S_\t$. This completes
the proof. \hfill $\square$

\medskip

\subsection{Proof of Theorem \ref{teo1}}

Before proving Theorem \ref{teo1}, we need some preliminary
results.

\begin{proposition}\label{prop5}
Let $\tilde w\in C^2([0,\infty))\cap L^\infty((0,\infty)).$ Let
\[w(x):=\tilde w(|x|)\quad \textrm{for all}\;\; x\in \re^N\,.\]
Set $r\equiv |x|$.If
\begin{equation}\label{e33} \tilde w''(r)+
\frac{N-2s+1}{r}\tilde w'(r)\geq 0,
\end{equation}
then $w$ is a supersolution to equation
\begin{equation}\label{e31}
(-\Delta)^s w\,=\,0\quad \textrm{in}\;\; \re^N\,.
\end{equation}
\end{proposition}
Observe that in Proposition \ref{prop5} $w$ is a supersolution to
equation \eqref{e31} in the sense of Definition \ref{defsole} with
$c\equiv 0$\,.

\noindent{\it Proof\,.} From \cite[Theorems 1.1, 1.2, and remarks
after Theorem 1.2]{FerrVerb}, due to \eqref{e33} we have:
\[-(-\Delta)^s w\geq 0\quad \textrm{in}\;\; \mathcal D'(\re^N),\]
$i.e.$,
\begin{equation}\label{e32}
\int_{\re^N}w(-\Delta)^s \zeta dx \leq 0
\end{equation}
for any $\zeta \in C^\infty_c(\re^N), \zeta\geq 0.$ Since $w\in
\mathcal L^s(\re^N)\cap C^{2s+\g}_{loc}(\re^N)$ for some $\g>0$,
from \eqref{e32} we can infer that
\begin{equation}\label{e34}
\int_{\re^N}\zeta (-\Delta)^s w dx \leq 0
\end{equation}
for any $\zeta \in C^\infty_c(\re^N), \zeta\geq 0$. Inequality
\eqref{e34} immediately yields the thesis. \hfill $\square$

\bigskip

In the sequel we shall use the next well-known result, concerning the
hypergeometric function $_{2}F_1(a,b,c,s)\equiv F(a,b,c,s)$, with
$a,b\in \re, c>0, s\in \re\setminus\{1\}$  (see \cite[Chapters 15.2, 15.4]{DLMF}).

\begin{lemma}\label{lemma2}
The following limits hold true:
\begin{itemize}
\item[(i)] if $c>a+b$, then
\[\lim_{s\to 1^-} F(a,b,c, s)=\frac{\Gamma(c)\Gamma(c-a-b)}{\Gamma(c-a)\Gamma(c-b)}\,;\]
\item[(ii)] if $c=a+b$, then
\[\lim_{s\to 1^-}\frac{F(a,b,c,s)}{-\log(1-s)}=\frac{\Gamma(a+b)}{\Gamma(a)\Gamma(b)}\,;\]
\item[(iii)] if $c<a+b$, then
\[\lim_{s\to 1^-}\frac{F(a,b,c,s)}{(1-s)^{c-a-b}}=\frac{\Gamma(c)\Gamma(a+b-c)}{\Gamma(a)\Gamma(b)}\,.\]
\end{itemize}
\end{lemma}
For further references, observe that
\begin{equation}\label{segnoga}
\Gamma(t)>0\quad\textrm{for all}\;\; t>0, \quad \Gamma(t)<0\quad
\textrm{for all}\;\; t\in (-1,0)\,.
\end{equation}
\bigskip

\noindent{\it Proof of Theorem \ref{teo1}\,.} Let $\psi=\psi(|x|)$
be defined as in \eqref{e66}, where $\beta>0$ is a constant to be
chosen. Also, let~$\alpha$ be as in~$(H_0)-(ii)$.
Set $r\equiv |x|$. We have:
\begin{equation}\label{e36}
\psi'(r)=-\b r(1+r^2)^{-\left(\frac{\b}2+1\right)}\quad
\textrm{for all}\;\; r>0\,,
\end{equation}
\begin{equation}\label{e37}
\psi''(r)=\b (1+r^2)^{-\left(\frac{\b}{2}+2\right)}[-1 +
(\b+1)r^2]\quad \textrm{for all}\;\; r>0\,.
\end{equation}
For any~$\lambda>0$ define
\[\phi(x,t):= e^{-\l t}\psi(r)\quad \textrm{for all}\;\; (x,t)\in \bar S_T\,.\]
At first observe that \eqref{e5} is satisfied.
\medskip

\noindent Suppose that $(i)$ is satisfied. In view of
\eqref{e36}-\eqref{e37}, we have:
\begin{equation}\label{e38}
\begin{split}
&\psi''(r)+\frac{N-2s
+1}{r}\psi'(r)\\
&\quad=\b(1+r^2)^{-\left(\frac{\b}2+2\right)}[(\b-N+2s)r^2-(N-2s+2)]\quad
\textrm{for all}\;\; r>0, t>0\,.
\end{split}
\end{equation}
Since $0<\b\leq N-2s$, by \eqref{e38},
\begin{equation}\label{e39}
\psi''(r)+\frac{N-2s +1}{r}\psi'(r)\leq 0 \quad \textrm{for
all}\;\; r>0\,.
\end{equation}
By Proposition \ref{prop5},
\begin{equation}\label{e40}
-(-\Delta)^s \phi(x,t)=- e^{-\l t}(-\Delta)^s\psi(r)\leq 0\quad
\textrm{for all}\;\; (x,t)\in \bar S_T\,.
\end{equation}
From \eqref{e40}, for any $\a\in \re$, we obtain
\begin{equation}\label{e41}
-(-\Delta)^s\phi(x,t) +\rho(x)\pa_t \phi(x,t)\leq -\l
\rho(x)e^{-\l t} < 0 \quad \textrm{for all}\;\; (x,t)\in \bar
S_T\,.
\end{equation}
By \eqref{e41} and Proposition \ref{prop1}, the conclusion
follows.

\medskip

In order to obtain the thesis of
Theorem \ref{teo1}
in cases $(ii),(iii),(iv)$, note
that (see the proof of Corollary 4.1 in \cite{FerrVerb}) we have:
\begin{equation}\label{e42}
-(-\Delta)^s\psi(r)= - \check C F(a,b,c, - r^2)\quad \textrm{for
all}\;\; r>1\,,
\end{equation}
where $\check C>0$ is a positive constant, and
\[a=\frac N 2 +s,\quad b= \frac{\b}2+s, \quad c=\frac N2\,.\]
By Pfaff's transformation,
\begin{equation}\label{e43}
F(a,b,c, -r^2)=\frac 1{(1+r^2)^b}F\left(c-a, b, c,
\frac{r^2}{1+r^2}\right)\quad \textrm{for all}\,\;r>1\,.
\end{equation}

\medskip

Suppose that $(ii)$ is satisfied. From Lemma \ref{lemma2}-$(i)$, \eqref{e42}
and \eqref{e43}, for any $\epsilon>0$, for some $R_\epsilon>1$, we
have:
\begin{equation}\label{e44}
-(-\Delta)^s\psi(r)\leq \check C(C_1 +
\epsilon)(1+r^2)^{-\left(s+\frac{\b}2\right)}\quad
\textrm{whenever}\;\;|x|>R_\epsilon\,,
\end{equation}
where $$C_1=-\frac{\Gamma(\frac N
2)\Gamma\left(\frac{N-\b}2\right)}{\Gamma\left(\frac{N+s}{2}
\right)\Gamma\left(\frac{N-\b}{2}-s\right)}>0$$ (see
\eqref{segnoga}). From \eqref{e44}, since $\a\leq 2s,$ we obtain
for all $|x|>R_\epsilon, t\in [0,T]$
\begin{equation}\label{e45}
\begin{split}
-(-\Delta)^s \phi(x,t)+\rho(x)\, \pa_t \phi(x,t) \hspace{3 cm}\\
\leq e^{-\l t}\big[(C_1+\epsilon)\check
C(1+r^2)^{-\left(s+\frac{\b}2\right)}-\l K
(1+r^2)^{-\frac{\a+\b}2}\big]<0\,,
\end{split}
\end{equation}
provided
\begin{equation}\label{e67}
\l>(C_1+\epsilon)\frac{\check C}K\,.
\end{equation}
On the other hand, for all $|x|\leq R_\epsilon,\; t\in [0,T]$
\begin{equation}\label{e48}
-(-\Delta)^s \phi(x,t) + \rho(x)\, \pa_t \phi(x,t) \hspace{3 cm}\\
\leq e^{-\l t}\big[ M_{\epsilon, \b} - \l K
(1+R^2_\epsilon)^{-\frac{\a+\b}2}\big] <0\,,
\end{equation}
taking
\begin{equation}\label{e68}
 \l>\frac{M_{\epsilon,\b}
(1+R^2_\epsilon)^{\frac{\a+\b}2} }{K},
\end{equation} where
$$M_{\epsilon,\b}:= \max_{x\in
\bar
B_{R_\epsilon}}\big\{\big|-(-\Delta)^s\psi(|x|)\big|\big\}\,.$$ By
\eqref{e45}, \eqref{e48} the conclusion follows by Proposition
\ref{prop1}.

\medskip

Suppose that $(iii)$ is satisfied. Let $\a<2s\,.$ From Lemma
\ref{lemma2}-$(ii)$ and \eqref{e43}, for any $\epsilon>0$, for
some $R_\epsilon>1$, we have:
\begin{equation}\label{e46}
-(-\Delta)^s\psi(r)\leq \check C(C_2 +
\epsilon)(1+r^2)^{-\left(s+\frac{\b}2\right)}\log(1+r^2)\quad
\textrm{whenever}\;\;|x|>R_\epsilon\,,
\end{equation}
where
$$C_2=-\frac{\Gamma\left(\frac{\b}2\right)}{\Gamma(-s)\Gamma\left(\frac{\b}2+s\right)}>0$$
(see \eqref{segnoga}). From \eqref{e46} we obtain for all
$|x|>R_\epsilon, t\in [0,T]$
\begin{equation*}
\begin{split}
-(-\Delta)^s \phi(x,t)+\rho(x)\, \pa_t \phi(x,t) \hspace{3 cm}\\
\leq e^{-\l t}\big[(C_2+\epsilon)\check
C(1+r^2)^{-\left(s+\frac{\b}2\right)}\log(1+r^2)-\l
(1+r^2)^{-\frac{\a+\b}2}\big]<0\,,
\end{split}
\end{equation*}
taking a possibly larger $R_\epsilon>1$, and so the desired claim
follows from \eqref{e48}.

Now, let $\a=2s$. From \eqref{e46} again, we have for all
$|x|>R_\epsilon, t\in [0,T]$
\begin{equation}\label{e47}
\begin{split}
-(-\Delta)^s \phi(x,t)+\rho(x)\, \pa_t \phi(x,t) \hspace{3 cm}\\
\leq e^{-\l t}\big[(C_2+\epsilon)\check
C(1+r^2)^{-\left(s+\frac{\b}2\right)}\log(1+r^2)-\l
(1+r^2)^{-\left(s+\frac{\b}2\right)}\log(1+r^2)\big]<0\,,
\end{split}
\end{equation}
provided
\begin{equation}\label{e69}
\l>(C_2+\epsilon)\frac{\check C}K\,.
\end{equation}

Furthermore, \eqref{e48} holds true, provided \eqref{e68} holds
true when $\a<2s$, or
\begin{equation}\label{e70}
\l>\frac 1{K s e}
\end{equation}
when $\a=2s$. From \eqref{e47} and \eqref{e48} the conclusion
follows.

\medskip

Finally, suppose that $(iv)$ is satisfied. From Lemma
\ref{lemma2}-$(iii)$ and \eqref{e43}, for any $\epsilon>0$, for
some $R_\epsilon>1$, we have:
\begin{equation}\label{e49}
-(-\Delta)^2\psi(r)\leq \check C(C_3 +
\epsilon)(1+r^2)^{-\left(s+\frac N 2\right)}\quad
\textrm{whenever}\;\;|x|>R_\epsilon\,,
\end{equation}
where $$C_3=-\frac{\Gamma(\frac N
2)\Gamma\left(\frac{\b-N}2\right)}{\Gamma(-s)\Gamma\left(\frac{\b}{2}+s\right)}>0$$
(see \eqref{segnoga}). From \eqref{e49}, since $\a+\b\leq 2s +N$,
we obtain for all $|x|>R_\epsilon, t\in [0,T]$
\begin{equation}\label{e50}
\begin{split}
-(-\Delta)^s \phi(x,t)+\rho(x)\, \pa_t \phi(x,t) \hspace{3 cm}\\
\leq e^{-\l t}\big[(C_3+\epsilon)\check C(1+r^2)^{-\left(s+\frac
N2\right)}-\l K (1+r^2)^{-\frac{\a+\b}2}\big]<0\,,
\end{split}
\end{equation}
provided
\begin{equation}\label{e71}
\l>(C_3+\epsilon)\frac{\check C}K\,.
\end{equation}
On the other hand, \eqref{e48} holds true, provided \eqref{e68} is
satisfied. In view of \eqref{e48} and \eqref{e50}, the conclusion
follows by Proposition \ref{prop1}. This completes the proof.
\hfill $\square$

\bigskip

\section{Elliptic equations: proofs}\label{pe}\setcounter{equation}{0}
\subsection{Proof of Proposition \ref{prop3}}
Analogously to Lemma \ref{lemma1} the next lemma can be shown.
\begin{lemma}\label{lemma1b}
Let $\phi\in C^2(\re^N), \phi>0$; suppose that \eqref{e7} is
satisfied. Let $v\in L^1_{\psi}(\re^N)$. Then
\[ \int_{\re^N}|v(x)|\phi(x)|(-\Delta)^s\g_R(x)|dx dt + \int_{\re^N}|v(x)|\,|\mathcal B(\phi, \g_R)(x)|dx \to 0 \]
as $R\to \infty\,.$
\end{lemma}

\noindent{\it Proof of Proposition \ref{prop3}}\,. Take a function
$v\in C^2(\re^N)$ with $supp\;\, v$ compact. Moreover, take a
function $w\in C(\re^N)$ such that $w\in \mathcal L^s(\re^N)\cap
C^{2s+\g}(\re^N)$ if $s<\frac 1 2$, or  $w\in \mathcal
L^s(\re^N)\cap C^{1, 2s+\g-1}_{loc}(\re^N)\; \textrm{if}\; s\geq
\frac 1 2$, for some $\g>0$. Integrating by parts we have:
\begin{equation}\label{e52}
\int_{\re^N} v \big[-(-\Delta)^s w -  \rho(x)c(x) w \big]\,  dx  =
\int_{\re^N} w\big[- (- \Delta)^s v-\rho(x)c(x) v\big]\, dx.
\end{equation}
Let $G_\a$ be defined as in \eqref{e23}. Thanks to \eqref{e1b} and
\eqref{e14} we obtain
\begin{equation}\label{e24}
\begin{split}
(-\Delta)^s[G_\a(u)] + \rho c G_\a(u) \leq p (u^2+\a)^{\frac p2
-1}u(-\Delta)^s u + \rho c G(u)\\
=p (u^2+\a)^{\frac p2 -1}u[(-\Delta)^s u + \rho c u] + \rho c
(u^2+\a)^{\frac p2}- p(u^2+\a)^{\frac p2 -1}u^2 c \rho\\
=(u^2+\a)^{\frac p2 -1}\rho c(u^2+\a-p u^2)=(u^2+\a)^{\frac
p2-1}\rho c[(1-p)u^2+\a]\quad \textrm{in}\;\; \re^N\,.
\end{split}
\end{equation}
From \eqref{e52} with $w=G_\a(u)$ and \eqref{e24} it follows that
\begin{equation}\label{e25}
\int_{\re^N} v (u^2+\a)^{\frac p2-1}\rho c[(p-1)u^2- \a]\leq
\int_{\re^N} G_\a(u)\big[- (- \Delta)^s v-\rho(x)c(x) v\big]\, dx.
\end{equation}
Letting $\a\to 0^+$ in \eqref{e25}, by the dominated convergence
theorem we get
\begin{equation}\label{e26}
\int_{\re^N}  |u|^p\big[(- \Delta)^s v +\rho(x) c(x) p v \big]
dx\leq 0\,.
\end{equation}
For any $R>0$, we can choose
\[v(x):= \zeta(x)\g_R(x)\quad \textrm{for all}\;\; x\in \re^N\,.\]
From \eqref{e9} we obtain
\begin{equation}\label{e27}
-(-\Delta)^s v- p \rho c v=\gamma_R \big[-(-\Delta)^s \zeta - p
\rho c \zeta\big]- \zeta(-\Delta)^s\g_R + \mathcal B(\zeta,\g_R)
\quad \textrm{in}\;\; \re^N\,.
\end{equation}
By \eqref{e26}, \eqref{e27},
\begin{equation}\label{e28}
\begin{split}
\int_{\re^N} |u(x)|^p \gamma_R \big[(-\Delta)^s \zeta + p \rho c
\zeta\big] dx \leq  \int_{\re^N}
|u|^p\big[ - \zeta(-\Delta)^s\g_R + \mathcal B(\zeta,\g_R) \big]dx \\
\leq \int_{\re^N}|u|^p\big[ |\zeta||(-\Delta)^s\g_R| + |\mathcal
B(\zeta,\g_R)|\big]dx\,.
\end{split}
\end{equation}
Hence, from Lemma \ref{lemma1b} and the monotone convergence
theorem, sending $R\to \infty$ in \eqref{e28} we get
\begin{equation}\label{e29}
\int_{\re^N} |u(x)|^p \big[(-\Delta)^s \zeta + p \rho c \zeta\big]
dx\leq 0\,.
\end{equation}
From \eqref{e29} and \eqref{e6}, since $|u|^p\geq 0$, we can infer
that $u\equiv 0$ in $\re^N$. This completes the proof.~\hfill
$\square$

\medskip

\subsection{Proof of Theorem \ref{teo3}}

\noindent{\it Proof of Theorem \ref{teo3}\,.} Let $\psi$ be
defined by \eqref{e66}. From the same arguments as in the proof of
Theorem \ref{teo1} we can infer that $\psi$ solves \eqref{e14},
for properly chosen $\b>0$. Note that to do this, we require that
by hypothesis $p c_0$ satisfies the same conditions as $\l$ in the
proof of Theorem \ref{teo3}. To be specific, we require that, for
some $\epsilon>0$,
\begin{equation}\label{e72}
p c_0 K >\max\big\{(C_1+\epsilon)\check C, M_{\epsilon,
\beta}(1+R^2_\epsilon)^{\frac{\a+\b}2}\big\}\,
\end{equation}
when $(ii)$ is satisfied (see \eqref{e67}, \eqref{e68});
\begin{equation}\label{e73}
p c_0 K >\max\big\{(C_2+\epsilon)\check C, M_{\epsilon,
\beta}(1+R^2_\epsilon)^{\frac{\a+\b}2}\big\}\,
\end{equation}
when $(iii)$ is satisfied and $\a<2s$ (see \eqref{e69},
\eqref{e68}), while
\begin{equation}\label{e74} p c_0 K
>\max\left\{(C_2+\epsilon)\check C, \frac 1{s e}\right\}\,
\end{equation}
if $\a=2s$ (see \eqref{e69}, \eqref{e70});
\begin{equation}\label{e75}
p c_0 K >\max\big\{(C_3+\epsilon)\check C, M_{\epsilon,
\beta}(1+R^2_\epsilon)^{\frac{\a+\b}2}\big\}\,\,.
\end{equation}
when $(iv)$ is satisfied (see \eqref{e71}, \eqref{e68}).
Thus, by Proposition \ref{prop3} the conclusion follows. \hfill
$\square$

\subsection{Proof of Theorem \ref{teo4}} In
order to prove Theorem \ref{teo4} we need the next lemma.
\begin{lemma}\label{l}
Let assumption $(H_0)-(i), (ii)$ be satisfied with $\alpha<s$. Let $N>-2s +\a.$
Then there exist constants $C>0$ and $\s>0$ such that, for any
$R>1$,
\begin{equation}\label{e45c}
R^\sigma \Big[\big|\phi(x)(-\Delta)^s \gamma_R(x)\big| +
\big|\mathcal B (\phi, \g_R)(x)\big|\Big]\leq  C \rho(x)\phi(x)
\quad \textrm{for all}\;\; x\in \re^N\,;
\end{equation}
here $\phi$ is defined as in \eqref{e76c}.
\end{lemma}
For any $0<r<R$ we set $B_R:=\{x\in \re^N\,:\, |x|<R\},\;
B_R^c:=\re^N\setminus B_R,\; \mathcal A_{r,R}:=B_R\setminus \bar
B_r$.

\medskip
\noindent{\it Proof\,.} Now we estimate $|\mathcal B(\phi,
\g_R)|$. To do this, we cover $\re^N\times\re^N$ as in the proof of Lemma \ref{lemma1}.
For any $R>1$ set
\[f_R(x,y):=\frac{|\phi(x,t)-\phi(y,t)|\,|
\g_{R}(x)-\g_{R}(y)|}{|x-y|^{N+2s}}\quad \textrm{for all}\;\;
x,y\in \re^N, x\neq y \,.\] Due to \eqref{e81}-\eqref{e82} we have
that for all $x\in B_{\frac R 2}$
\begin{equation}\label{e8c}
\int_{B_{\frac R 2}} f_R(x,y) dy\,=\,0;
\end{equation}
furthermore, for all $x\in B_R^c$
\begin{equation}\label{e9c}
\int_{B_R^c} f_R(x,y) dy\,=\,0.
\end{equation}

\smallskip
Note that in the sequel we shall denote by the same $C$ different
positive constants independent of $R.$ Let $0<\beta<N,
0<\sigma<\min\{\beta, 2s\}$ be two parameters to be chosen later. For all
$(x,y)\in A_1$ we get
\begin{equation}\label{e10c}
|x-y|\geq C |x|\,,
\end{equation}
and
\begin{equation}\label{e11c}
|x-y|\geq\frac{R}{2}-\frac{R}{8} \geq  \frac R 4+|y|\,.
\end{equation}
Hence, \begin{equation}\label{e12c}
\frac 1{|x-y|^{N+2s}}\leq \frac C{|x|^{\b}\big( \frac R 4+
|y|\big)^{N+2s-\beta}}\,.
\end{equation}
Moreover, from \eqref{e7c} it follows that, if~$(x,y)\in A_1$,
then, for any $R>1$,
\begin{equation}\label{e13c}
|\phi(x)|+|\phi(y)|\leq \frac{ C}{1+|y|^{N-2s}}.
\end{equation}
Inequalities \eqref{e12c} and \eqref{e13c} yield, for any $x\in
B_{\frac R 2}^c$,
\begin{equation}\label{e14c}
R^\s\int_{B_{\frac R 8}} f_R(x,y) dy \leq C
\frac{R^\s}{|x|^\b}\int_{B_{\frac R
8}}\frac{dy}{1+|y|^{2N-\b}}\leq \frac C{1+|x|^{\b-\s}}\,.
\end{equation}

\smallskip

For all $(x,y)\in A_2$ we have
\begin{equation}\label{e15c}
|\phi(x)|+|\phi(y)|\leq \frac{ C}{1+|x|^{N-2s}},
\end{equation}
and
\begin{equation}\label{e16c}
|x-y| \geq C|y|\,.
\end{equation}
In view of \eqref{e15c} and \eqref{e16c}, we obtain for any $x\in
B_{\frac R 8}$
\begin{equation}\label{e16bc}
R^\s \int_{B^c_{\frac R 2}} f_R(x,y) dy \leq \frac{C
R^\s}{1+|x|^{N-2s}}\int_{B^c_{\frac R 2}}\frac{dy}{|y|^{N+2s}}
\leq \frac{C R^{-2s +\s}}{1+|x|^{N-2s}}\leq \frac
C{1+|x|^{N-\s}}\,.
\end{equation}

\smallskip
Also, for all $(x,y)\in A_3$ we have that \eqref{e10c},
\eqref{e13c} and \eqref{e16c} hold true. From \eqref{e10c} and
\eqref{e16c} we get
\begin{equation}\label{e17c}
\frac 1{|x-y|^{N+2s}}\leq \frac C{|x|^\b|y|^{N+2s-\b}}\,.
\end{equation}
So, due to \eqref{e13c} and \eqref{e17c}, we obtain for any $x\in
B_{2R}^c$
\begin{equation}\label{e18c}
R^\s\int_{\mathcal A_{\frac R 8, R}}f_R(x,y) dy \leq \frac
C{|x|^\b}\int_{\mathcal A_{\frac R 8,
R}}\frac{R^\s}{|y|^{2N-\b}}dy\leq \frac C{1+|x|^{\b-\s}}\,.
\end{equation}
For all $(x,y)\in A_4$, we have that  \eqref{e15c} and
\eqref{e16c} hold true. Hence, for any $x\in \mathcal A_{\frac R
8, R}$,
\begin{equation}\label{e19c}
R^\sigma\int_{B^c_{2R}}f_R(x,y) dy\leq \frac
C{1+|x|^{N-2s}}\int_{B^c_{2R}}\frac{R^\s}{|y|^{N+2s}}dy\leq
\frac{C R^{-2s+\s}}{1+|x|^{N-2s}}\leq \frac{C}{1+|x|^{N-\s}}\,.
\end{equation}

\medskip

Now, let $(x,y)\in A_5$. We shall distinguish the cases $s\in \left(0,\frac 12\right]$ and $s\in\left(\frac 12 ,1\right)$. To begin with, take any $s\in\left(0,\frac 12\right]$. Since
in $A_5$ the roles of $x$ and $y$ are symmetric, from \eqref{e7c}
we can infer that
\begin{equation}\label{e20c}
|\phi(x)|+|\phi(y)|\leq \frac{C}{|x|^{N-2s}}.
\end{equation}
Furthermore,
\begin{equation}\label{e21c}
|\g_{R}(x)-\g_{R}(y)|\leq \frac{C}{R}|x-y|.
\end{equation}
Thus, for any $0<\d<\frac{N-\beta}{N+1-2s}$, for all $(x,y)\in A_5$,
\[|\phi(x)-\phi(y)|=|\phi(x)-\phi(y)|^\d |\phi(x)-\phi(y)|^{1-\d}\leq \frac{C |x-y|^\d}{|x|^{(N-2s)(1-\d)}}\,.
\]

Hence, for any $x\in \mathcal A_{\frac R 8, R}$,
\begin{equation}\label{e22c}
R^\s\int_{\mathcal A_{\frac R 8, R}}f_R(x,y) dy\leq \frac{C
R^{\s-1}}{|x|^{(N-2s)(1-\d)}}\int_{\mathcal A_{\frac R 8, R}}\frac
1{|x-y|^{N+2s-1-\d}}dy\,.
\end{equation}
By the change of variables $\tilde y:=x-y$, since $s\in \left(0,
\frac 12\right]$, from \eqref{e22c} it follows that, for all $x\in \mathcal A_{\frac R8, R}$,
\begin{equation}\label{e23c}
\begin{split}
R^\s\int_{\mathcal A_{\frac R 8, R}}f_R(x,y) dy\leq \frac{C
R^{\s-1}}{|x|^{(N-2s)(1-\d)}}\int_{\frac R 8\leq |x-\tilde y|\leq 2R}\frac
1{|\tilde y|^{N+2s-1-\d}}d\tilde y\\
\leq \frac{C R^{\s-1}}{|x|^{(N-2s)(1-\d)}}\int_{B_{4R}}\frac{1}{|\tilde
y|^{N+2s-1-\d}}d\tilde y\leq C
\frac{R^{-2s+\s+\d}}{|x|^{(N-2s)(1-\d)}}\\ \leq \frac C{1+ |x|^{N-\s-\d(N+1-2s)}}\leq 
 \frac C{1+ |x|^{\b-\s}}
\,.
\end{split}
\end{equation}
Now, let $s\in\left(\frac 12 ,1\right)$. By \eqref{e7c}, we get
\begin{equation}\label{e28c}
|\phi(x)-\phi(y)|\leq\frac{C}{1+|z|^{N-2s}}|x-y|,
\end{equation}
for some $z$  in the segment joining $x$ and $y$. For any $R>0$
let
$$Q_R\equiv Q:=\left\{(x,y)\in A_5:\, |x-y|\leq \frac{R}{100}\right\}.$$
Note that, if $(x,y)\in Q$ then every point $z$ lying on the
segment from $x$ to $y$ verifies $|z|\geq C|x|$. Hence, since
$s\in \left( \frac 1 2,1\right)$, \eqref{e28c} and \eqref{e21c}
yield, for all $x\in \mathcal A_{\frac R 8, 2R}$,
\begin{equation}\label{e29c}
\begin{split}
R^\s\int_{\frac R 8<|y|<2R, |x-y|<\frac R{100}}f_R(x,y) dy\leq
\int_{\frac R 8<|y|<2R, |x-y|<\frac R{100}}\frac{C
R^{\s-1}}{|x|^{N-2s}|x-y|^{N+2s-2}} dy\\
\leq \frac{C R^{\s-1}}{|x|^{N-2s}}\int_{B_{4R}}\frac{d\tilde
y}{|\tilde y|^{N+2s-2}}\leq \frac{C R^{1-2s+\s}}{|x|^{N-2s}}\leq
\frac{C}{1+|x|^{N-\s}}\,.
\end{split}
\end{equation}
On the other hand, if 
$(x,y)\in A_5\setminus Q$ we have that
\begin{equation}\label{e30c}
|x-y|>\frac{R}{100}\geq C|y|.
\end{equation}
Then, by \eqref{e20c} and \eqref{e30c}, for all $x\in \mathcal
A_{\frac R 8, 2R}$,
\begin{equation}\label{e31c}
\begin{split}
R^\s\int_{\frac R 8<|y|<2R, |x-y|>\frac R{100}}f_R(x,y) dy\leq
\frac{C R^{\s-1}}{|x|^{N-2s}}\int_{\mathcal A_{\frac R 8,
R}}\frac{dy}{|y|^{N+2s-1}}\\ \leq \frac{C
R^{-2s+\s}}{|x|^{N-2s}}\leq \frac{C}{1+|x|^{N+1-\s}}\,.
\end{split}
\end{equation}
Therefore, from \eqref{e23c}, \eqref{e29c} and \eqref{e31c}, we have for all
$x\in \mathcal A_{\frac R 8, 2R} $, for each $s\in (0, 1)$,
\begin{equation}\label{e32c}
R^\sigma\int_{\mathcal A_{\frac R 8, 2R}} f_R(x,y) dy \leq \frac
C{1+|x|^{N-\s}}\,.
\end{equation}
From \eqref{e8c}, \eqref{e9c}, \eqref{e14c}, \eqref{e16bc}, \eqref{e18c}, \eqref{e19c} and \eqref{e32c} we obtain, for all $x\in \re^N$,
\begin{equation}\label{e33c}
R^\s\int_{\re^N} f_R(x,y) dy\leq \frac{C}{1+|x|^{\b-\s}}\,.
\end{equation}
Since $N-2s+2\a>0$, we can choose $0<\beta<N-2s+2\a$, and then
\begin{equation}\label{e33bc}
0<\sigma<\beta-N+2s-2\alpha.
\end{equation}
Thus
\begin{equation}\label{e34c}
\b-\s>2\a+N-2s\,.
\end{equation}
Note that by $(H_0)-(ii)$ and \eqref{e7c},
\begin{equation}\label{e35c}
\rho(x)\phi(x)\geq \frac{C}{1+|x|^{2\a+N-2s}}\quad \textrm{for
all}\;\; x\in \re^N\,.
\end{equation}
In view of \eqref{e33c}, \eqref{e34c} and \eqref{e35c} we can
infer that
\begin{equation}\label{e36c}
R^\s\big|\mathcal B(\phi, \g_R)(x) \big|\leq C \rho(x)\phi(x)\quad
\textrm{for all}\,\, x\in \re^N\,,
\end{equation}
with $\s$ as in \eqref{e33bc}.

\medskip
Observe that, for all $x\in \re^N$,
\begin{equation}\label{e37c}
\big|(-\Delta)^{s}\gamma_{R}(x)\big|=R^{-2s }\left|
(-\Delta)^{s}\gamma \left(\frac{x}{R}\right)\right|\leq
C R^{-2s}.
\end{equation}
Take $\sigma>0$ satisfying \eqref{e33bc}. For any $R>1, x\in B_{R}$ we get
\begin{equation}\label{e38c}
R^\sigma\big|(-\Delta)^s\g_R(x)\big|\leq \frac{C R^\s}{R^{2s}}\leq
\frac{C}{1+R^{2s-\s}}\leq
\frac{C}{1+|x|^{2s-\s}}\,.
\end{equation}
Note that for all $x\in \re^N$
\[|(-\Delta)^s \g_1(x)|\leq \frac C{1+ |x|^{N+2s}},\]
so
\begin{equation}\label{e39c}
|(-\Delta)^s \g_R(x)|\leq \frac {C R^{-2s}}{1+
\left(\frac{|x|}{R}\right)^{N+2s}}\,.
\end{equation}
This implies that for any $x\in B^c_{R}$ we have
\begin{equation}\label{e40c}
R^\s|(-\Delta)^s\g_R(x)|\leq \frac{C R^{N+2s+\s}}{R^{2s}(R^{N+2s}+|x|^{N+2s})}\leq
\frac{C}{1+|x|^{2s-\s}}\,.
\end{equation}
Choose $\s>0$ so that
\begin{equation}\label{e41c}
0<\s<2s-2\a\,.
\end{equation}
Thus,
\begin{equation}\label{e42c}
2\a<2s-\s\,.
\end{equation}
In view of \eqref{e38c},
\eqref{e40c} and \eqref{e42c} we obtain
\begin{equation}\label{e43c}
R^\s|(-\Delta)^s\g_R(x)|\leq \frac C{1+|x|^{2\a}}\quad \textrm{for
all}\;\; x\in \re^N\,.
\end{equation}
By \eqref{e43c} and $(H_0)-(ii)$,
\begin{equation}\label{e44c}
R^\s\phi(x)|(-\Delta)^s\g_R(x)|\leq C \phi(x) \rho(x)\quad
\textrm{for all}\;\; x\in \re^N\,.
\end{equation}
We can select $\sigma>0$ such that both \eqref{e33bc} and
\eqref{e41c} hold true. From \eqref{e36c} and \eqref{e44c} we get
\eqref{e45c}. This completes the proof. \hfill $\square$

\bigskip
\medskip

Let us prove Theorem \ref{teo4}.

\noindent{\it Proof of Theorem \ref{teo4}} Let $u$ be any solution
to equation \eqref{e1b}. We can repeat the proof of Proposition
\ref{prop3} in order to obtain \eqref{e28}. Then we choose
$\zeta=\phi$. This combined with \eqref{e6c} yields
\begin{equation}\label{e48c}
\begin{split}
\int_{\re^N} |u(x)|^p\g_R(x)F(x) dx + p \int_{\re^N}|u(x)|^p
\phi(x)\g_R(x)c(x)\rho(x) \,dx \\ \leq  - \int_{\re^N}
|u(x)|^p\Big\{\phi(x)(-\Delta)^{s}[\g_R(x)] - \mathcal
B(\phi,\g_R)(x) \Big\}\, dx=:\mathcal I(R)\,.
\end{split}
\end{equation}
Since $ \g_R\geq 0, \phi>0, |u|^p\geq 0,
\rho>0$, due to \eqref{e48c} we obtain
\begin{equation}\label{e49c}
\int_{\re^N} |u(x)|^p\g_R(x)F(x)dx dt \\ \leq  \mathcal I(R)\quad
\textrm{for all}\;\; R>1\,.
\end{equation}

\medskip

We claim that
\begin{equation}\label{e71c}
\liminf_{R\to \infty} \mathcal I(R)\,=\,0\,.
\end{equation}

\smallskip

In fact, suppose by contradiction that $\kappa:=\liminf_{R\to
\infty} \mathcal I(R)>0$. So, there exists $R_0>0$ such that
\begin{equation}\label{e73c}
\mathcal I(R) \geq \frac{\kappa}2\quad \textrm{for all}\;\;
R>R_0\,.
\end{equation}
From \eqref{e45c}, \eqref{e8} and \eqref{e73c} we have, for all
$R>R_0$,
\begin{equation}\label{e75c}
\begin{split}
p\int_{\re^N} |u(x)|^p c(x) \rho(x) \phi(x) dx  \\ \geq
\frac{c_0}{C} R^{\s}\int_{\re^N}|u(x)|^p
\Big[\big|\phi(x)(-\Delta)^s
\gamma_R(x)\big| + \big|\mathcal B (\phi, \g_R)(x)\big|\Big]dx \\
\geq \frac{c_0}{C} R^{\s} \mathcal I(R) \geq \frac{c_0
\kappa}{2C}R^{\s}\,.
\end{split}
\end{equation}
Sending $R\to \infty$ in \eqref{e75c}, we deduce that
$\int_{\re^N} c(x) \rho(x) |u(x)|^p \phi(x) dx = \infty\,.$ This
is in contrast with the hypothesis $u\in L^p_{c\rho\phi}(\re^N).$
Thus, $\kappa=0$, and the Claim is proved.

\medskip

From \eqref{e49c} and \eqref{e71c}, by Fatou's Lemma,
\[p\int_{\re^N} |u(x)|^p c(x) \rho(x) F(x) dx \,\leq \,0\,.\]
Since $|u|^p\geq 0, c>0, \rho>0$, and $F\geq 0$ were arbitrary we
can infer that $u\equiv 0$ in $\re^N\,.$ This completes the proof.
\hfill $\square$


\begin{thebibliography}{999}

\bibitem{AK}  H. Abels, M. Kassmann, {\it The Cauchy problem and the martingale problem for integro-differential operators with non-smooth kernels}, Osaka J. Math., {\bf 46 } (2009)\,, 661--683\,.

\bibitem{AB} D.G. Aronson, P. Besala, {\it Uniqueness of solutions to the Cauchy problem for parabolic equations}, J. Math.
Anal. Appl. {\bf 13} (1966), 516--526\,.

\bibitem{BPSV} B. Barrios, I. Peral, F. Soria, E. Valdinoci, {\it A Widder's type theorem for the heat equation with nonlocal
diffusion}, preprint (2013),
available on-line at
http://arxiv.org/abs/1302.1786\,.

\bibitem{BG} R.M. Blumenthal, R.K. Getoor, {\it Some theorems on stable processes}, Trans. Amer. Math. Soc. {\bf 95}  (1960), 263--273\,.

\bibitem{DPVal} E. Di Nezza, G. Palatucci, E. Valdinoci, {\it Hitchhiker's guide to the fractional Sobolev spaces}, Bull. Sci. Math. {\bf 136} (2012), 521--573\,.

\bibitem{EKP} S. D. Eidelman, S. Kamin and F. Porper, \emph{Uniqueness of solutions of the Cauchy problem for  parabolic equations degenerating at infinity},  Asympt. Anal., \textbf{22} (2000), 349--358\,.

\bibitem{FFV} F. Ferrari, B. Franchi, I.E. Verbitsky, {\it Hessian inequalities and the fractional
Laplacian}, J. Reine Angew. Math. (to appear)\,.

\bibitem{FerrVerb} F. Ferrari, I. E. Verbitsky, {\it Radial fractional Laplace operators and hessian
inequalities}, preprint (2012),
available on-line at
http://arxiv.org/abs/1203.3149\,.

\bibitem{IKO} A. M. Il'in, A. S. Kalashnikov, O. A. Oleinik, {\it Linear equations of the second order of parabolic
type}, Russian Math. Surveys {\bf 17} (1962), 1--144\,.

\bibitem{KPT} S. Kamin, M.A. Pozio, A. Tesei, {\it Admissible conditions for parabolic equations degenerating at
infinity}, St. Petersburg Math. J. {\bf 19} (2008), 239--251\,.

\bibitem{MP1} R. Mikulevicius, H. Pragarauskas, {\it On the Cauchy problem for integro-differential operators in Holder classes and the uniqueness of the martingale problem}\,, Potential Anal. (to appear) doi: 10.1007/s11118-013-9359-4.

\bibitem{MP2} R. Mikulevicius, H. Pragarauskas, {\it On the Cauchy problem for integro-differential operators in Sobolev classes and the martingale problem}\,, arXiv:1112.4467 (2011).


\bibitem{OR} O. A. Oleinik and E. V. Radkevic, ``Second Order Equations with
Nonnegative Characteristic Form", Amer. Math. Soc., Plenum Press,
New York - London, 1973\,.

\bibitem{DLMF} F.W.J. Olver, D. W. Lozier, R.F. Boisvert, C.W. Clark (eds.), "NIST Handbook of Mathematical functions", Cambridge University Press, New York, NY, 2010,
available on-line at http://dlmf.nist.gov\,.


\bibitem{P2} F. Punzo, {\it Uniqueness of solutions to degenerate parabolic and elliptic
equations in weighted Lebesgue spaces}, Math. Nachr. {\bf 286} (2013), 1043--1054\,.

\bibitem{PT} F. Punzo, G. Terrone, {\it On the Cauchy problem for a general fractional porous medium equation with variable density}, Nonlin. Anal. 
{\bf 8} (2014), 27--47\,.

\bibitem{Silv} L. Silvestre, ``Regularity of the obstacle problem for a fractional power of the Laplace operator" (PhD Thesis). The
University of Texas at Austin (2005)\,.

\bibitem{Ti} A. N. Tihonov, {\it
Th\'eor\`emes d'unicit\'e pour l'\'equation de la chaleur},
Mat. Sb. {\bf 42} (1935), 199--215\,.

\end{thebibliography}
\end{document}